\documentclass[11pt,final,
hypertex, 
titlepage]{article}

\usepackage{amsmath,verbatim,amssymb}
\usepackage{latexsym}
\usepackage{graphicx}
\usepackage{color}
\usepackage{multirow}
\usepackage{array}
\usepackage{booktabs}
\usepackage{subfigure}
\usepackage{url}
\usepackage{amsthm}
\usepackage{wrapfig}
\newcommand{\suchthat}{\mbox{{\bf $\ \mid \ $}}}

\newcommand{\abs}[1]{\lvert#1\rvert}

\newcommand{\ls}[1]
{\dimen0=\fontdimen6\the\font \lineskip=#1\dimen0
\advance\lineskip.5\fontdimen5\the\font \advance\lineskip-\dimen0
\lineskiplimit=.9\lineskip \baselineskip=\lineskip
\advance\baselineskip\dimen0 \normallineskip\lineskip
\normallineskiplimit\lineskiplimit \normalbaselineskip\baselineskip
\ignorespaces }

\newcommand{\Pf}{\paragraph{{\bf Proof.}}}       
\newcommand{\blot}{\hfill{\vrule height .9ex width .8ex depth -.1ex }}
\newcommand{\EndPf}{\hfill $\blot$ \medskip}     

\iftrue 
\usepackage{amsmath}
\usepackage{amsfonts}
\usepackage{verbatim}
\usepackage{tikz}

\usetikzlibrary{arrows,decorations.pathmorphing,backgrounds,positioning,fit}

\newcommand{\field}[1]{\mathbb{#1}}
\DeclareMathOperator{\PR}{\field{P}}             
\DeclareMathOperator{\E}{\field{E}}              
\def\N{\field{N}}                                
\def\R{\field{R}}                                
\def\F{\field{F}}                                
\def\A{\field{A}}                                
\def\X{\field{X}}                                

\def\Z{\field{Z}}
\else
\def\PR{\mathop{\rm I\kern -0.20em P}\nolimits}  
\def\E{\mathop{\rm I\kern -0.20em E}\nolimits}   
\def\N{\mathop{\rm I\kern -0.20em N}\nolimits}   
\def\R{\mathop{\rm I\kern -0.20em R}\nolimits}   
\def\F{\mathop{\rm I\kern -0.20em F}\nolimits}   
\def\A{{A}}                                      
\def\X{{X}}                                      
\fi

\newtheorem{thm}{Theorem}[section]

\newtheorem{prop}[thm]{Proposition}
\theoremstyle{definition}
\newtheorem{ex}{Example}[section]
\theoremstyle{plain}
\newtheorem{defn}[thm]{Definition}

\numberwithin{equation}{section}

\graphicspath{{./}{Figs/}}

\usepackage[colorlinks=true,linkcolor=magenta,citecolor=blue,pagebackref]{hyperref}

\title{Dynamic Service Rate Control for a Single Server Queue with Markov Modulated Arrivals}

\author{Ravi Kumar\thanks{rk454@cornell.edu}, Mark E. Lewis\thanks{mark.lewis@cornell.edu}
and Huseyin Topaloglu\thanks{ht88@cornell.edu}\\
{\small \it School of  Operations Research and Information Engineering}\\
{\small \it Cornell University}\\ {\small
\it Ithaca, NY 14853} \\
}
\date{\small \today}

\begin{document}
\maketitle
\setlength{\unitlength}{1mm}
\begin{abstract}
We consider the problem of service rate control of a single server queueing system with
a finite-state Markov-modulated Poisson arrival process. We show that the optimal
service rate is non-decreasing in the number of customers in the system; higher
congestion rates warrant higher service rates. On the contrary, however, we show that
the optimal service rate is not necessarily monotone in the current arrival rate. If the
modulating process satisfies a stochastic monotonicity property the monotonicity is
recovered.

We examine several heuristics and show where heuristics are reasonable substitutes for
the optimal control. None of the heuristics perform well in all the regimes. Secondly, we
discuss when the Markov-modulated Poisson process with service rate control can act as
a heuristic itself to approximate the control of a system with a periodic
non-homogeneous Poisson arrival process. Not only is the current model of interest in
the control of Internet or mobile networks with bursty traffic, but it is also useful in
providing a tractable alternative for the control of service centers with non-stationary
arrival rates.
\end{abstract}

\section{Introduction}
\label{sec:intro}

In this paper, we study a fundamental queueing control problem; managing the service
rate of a server in the face of non-stationary arrival rates. We have a queue with an
infinite buffer and a single server. Arrivals occur according to a Markov-modulated
Poisson process (MMPP), which is to say that the rate of the Poisson process driving the
arrivals into the system changes according to an exogenous Markov process. This
exogenous Markov process is commonly referred to as the \emph{phase modulating
process}. The service times are assumed to be exponential.\footnote{In Kendall's notation, our queueing
system is classified as MMPP/M/1}~We incur a holding cost for each job in the system
and there is a cost for running the server at different service rates.~Given that the state
of the phase modulating process and state of the queue are  known, we want to find a
policy to adjust the service rate so as to minimize either the expected discounted cost or
the long-run average cost rate over an infinite horizon.

Our model is motivated by the power aware transmission policies that are becoming
increasingly important over the Internet and in mobile networks. The goal of such policies
is to control the power consumption of wireless devices by adjusting the transmission rate
in response to the number of packets waiting to be transmitted in the buffer. Due to
changes in incoming and outgoing traffic through the device, it is almost always the case
that the packet arrivals display non-stationarities, creating periods of bursts followed by
near-complete silence.~Our use of an MMPP to model arrivals is intended to capture such
periods of bursts and silence. An alternative model to capture the non-stationary nature
of arrivals is the \mbox{non-homogenous} Poisson process (NHPP), but for the wireless
applications we have in mind, MMPP appears to be a more suitable model of
non-stationarity since these applications involve arrival rate changes occurring at random
points in time, whereas an NHPP models arrival rates as a fixed function of time.
Furthermore, when computing the optimal policy under an NHPP, one needs to keep track
of time together with the queue length, resulting in an uncountable state space. This
issue is not present when dealing with an MMPP.

Controlling queues when arrivals have varying rates poses interesting challenges. When
controlling such queues, the policy in use not only needs to consider the current arrival
rate, but it also needs to anticipate the arrival rate in the near future and adjust
decisions accordingly. For example, if the current arrival rate is relatively low, but arrivals
are expected to be more frequent in the near future, then the control policy may choose
to proactively speed up the service rate to empty the system (as much as possible)
before the higher arrival rate strikes. Similarly, if the current arrival rate is high and the
current system load is high, the control policy may slow down the service rate in
anticipation of lower arrival rates in the near future. The extent to which changes in
arrival rates can be foreseen or anticipated depends on how the non-stationarity is
modeled, but policies that explicitly address the non-stationarity in arrival rate are
naturally expected to make better decisions than those that do not.~Furthermore, the
need to model arrival non-stationarities is becoming increasingly important as
non-stationary queues find more use in telecommunications applications to study
congestion problems on the Internet and mobile networks.

We provide a characterization of how the optimal policy depends on the queue length and
the arrival rate. Throughout, we assume that the arrival intensities change according to a
general continuous time Markov chain (CTMC) defined on the state space
$\{1,2,\dots,L\}$ and the arrival rates of the MMPP are ordered such that
$\lambda_1\leq\lambda_2\leq\dots\leq\lambda_L$. In this context, the first interesting
question is whether the optimal service rate is monotone in the queue length, for a fixed
state of the phase modulating process. We answer this question, not too surprisingly, in
the affirmative, indicating that the optimal service rate is higher as we have more jobs in
the buffer, all else being equal. The second interesting question is whether the optimal
service rate is monotone in the state of the phase modulating process, for a fixed queue
length.~Perhaps surprisingly, the answer to this question is not necessarily affirmative,
indicating that the optimal service rate is not necessarily higher as the arrival rate
becomes larger, all else being equal. This observation builds the intuition that the optimal
policy  should anticipate the arrival rate in near future. For example, even if the current
arrival rate is higher, the optimal policy may choose not to serve the jobs faster because
the arrival rate is expected to slow down soon after the higher arrival rate strikes. Thus,
although it may be surprising to see that the optimal service rate is not necessarily
monotone in the state of the phase modulating process, this non-monotonicity embodies
the intuitive expectation that the optimal policy may start using faster service rates even
before higher arrival rates strike or may start using slower service rates even before
arrival rates slow down. Motivated by this observation, a natural question is when we can
expect the optimal policy to actually be monotone in the state of the phase modulating
process so that the behavior that we just mentioned is not prevalent. We give sufficient
conditions under which the optimal service rate is indeed monotone in the state of the
phase modulating process. These conditions are simple to check and they only depend on
the structure of the CTMC driving the phase process. These structural results are not only important in providing insights but are also useful in deriving efficient approximation methods. When the phase process for the MMPP has a large number of states, computing an optimal policy using value or policy iteration may still be a difficult task. In these situations, the structural properties of the value function can be used to develop approximate dynamic programming methods and obtain approximate results efficiently (see for example Powell\cite{Powell}).

We include a numerical study with two goals in mind. First, we examine when it is important to explicitly capture the non-stationary behavior of an arrival process via MMPP as opposed to using some natural heuristic like assuming the system has stationary arrivals. To implement the optimal policy, a decision maker needs to look at both the state of the queue and the state of the phase process of the MMPP while a heuristic control mechanism based only on the queue length or some fixed service rate may be easier to implement. Thus, a comparison between the two helps in determining the value of a more complex control mechanism. Second, since we mentioned the alternative of using an NHPP to capture the non-stationarity in the arrivals, we explore the possibility of using a ``suitable'' MMPP to approximate the control policy for a system with an NHPP with a periodic rate function. We find that our preliminary results are encouraging. This is a significant diversion from previous studies since the focus here is on computing an optimal control and not solely on evaluating performance measures.

Most of the research related to Markov-modulated queueing systems deals
with performance characteristics for systems without control. Excellent overviews of this line of work can be found in the survey
paper by Prabhu \cite{Prabhu1989} or more recently in Gupta et al.
\cite{Gupta2006}. A hierarchical scheme based on MMPP was proposed by
Muscariello et al. \cite{Muscariello2004} to model the data generated by Internet
users. Heffess et al. \cite{Heffes1986} used MMPP to approximate a statistical
multiplexer whose inputs consist of superposition of packetized voice sources and
data. A two state MMPP model was proposed by Shah-Heydari
\cite{Shah-Heydari1998} to model the so-called aggregate asynchronous transfer mode
traffic. For more general scenarios, Frost \cite{Frost1994} proposed a scheme to
approximate a simple NHPP using an MMPP by suitably quantizing the rate function of
the NHPP into a finite number of rates. Each rate corresponds to a state in the Markov
modulating process and the parameters of the MMPP model can be estimated using
empirical data.
%
%

There is also a rich body of literature on the subject of monotone optimal policies for
the control of a single server queue in a setting similar to the one considered here but
with stationary arrivals. See, for example, the classical work of Crabill \cite{Crabill1974},
Lippman \cite{Lippman1975} and Stidham and Weber \cite{StidhamJr1989}. In the
context of telecommunications systems, the existing literature addresses a more
closely related problem of service rate control of queues when the job service
requirements are influenced by an exogenous stochastic process\footnote{Similar to
the present work, the policy for this type of model depends on both the queue length
and the state of the exogenous process.}. Such models arise frequently in
point-to-point wireless data transmission where the induced transmission rates are
affected by the time varying properties of the transmission medium. Berry
\cite{Berry2000} considers a very general model for this problem under a
discrete-time Markovian setting. In this work, packet arrivals follow a batch Markov
process and the state of the transmission channel varies according to a secondary
discrete-time Markov chain. The data buffer and transmitter are modeled using a
single server queue with finite capacity. The goal of the transmitter is to minimize the
average cost rate or power consumption over an infinite time horizon subject to a
constraint on packet delay.~Another case with a constraint on the probability
of buffer overflow is also discussed in this work. The author proves several results related to
the monotonicity of the optimal policy. Motivated by mobile networks, Ata and
Zachariadis \cite{Ata2006} address the problem of finding optimal service rates for
multiple users that are being served by a central controller. Data gets transmitted
through a time varying channel that is being modulated by a two-state continuous
time Markov chain.~Packet data for each user arrives based on a Poisson process and
gets stored in a finite capacity queue before getting transmitted. The objective is to
maximize some measure of overall quality of service subject to a constraint on the
long-run average power consumption.~The authors show that the optimal service
rates for each user depend only on its own queue length and the state of the
transmission channel. They also present a method to explicitly characterize the
optimal policy for each user. To the best of our knowledge, none of the
aforementioned work considers the case of a multi-state Markov-modulated arrival
process with service rate control as is discussed in this paper.

We have organized this paper as follows. In Section
\ref{sec:formulate} we give a detailed description of the model and the associated
assumptions, and formulate the problem as a Markov decision process. In the average
cost case we provide stability conditions that guarantee convergence of the
discounted cost value function to the relative value function. In Section
\ref{sec:structure} we give structural results related to the optimal
policy in each case. 
In Section \ref{sec:numeric}, we present a detailed numerical study comparing the
performance of the optimal policy with heuristic policies and present an exploratory
study related to the computation of a heuristic policy for non-homogeneous Poisson
arrivals using the optimal policy for a MMPP/M/1 queue. We conclude the paper in
Section \ref{sec:conclude}.

\section{Model Formulation}\label{sec:formulate}

We consider a single server queue with infinite buffer capacity and job arrivals that
follow an MMPP.~Each arriving job has an
exponentially distributed service requirement with mean 1.~The phase transition
process for arrivals is an \textit{ergodic}, finite state continuous time Markov chain
with generator matrix $Q$. Let the state space for this process be denoted by
$\mathcal{S}:=\{1,2,\dots,L\}$. When the phase transition process is in phase
$s\in\mathcal{S}$, jobs arrive to the queue according to a Poisson process with rate
$\lambda_s$. Without loss of generality we assume that the states are ordered such
that $\lambda_1\leq\lambda_2\leq\dots\leq\lambda_L$. Let the number of jobs in the
system (buffer state) be denoted by $n \in \Z^+$, where $\Z^+$ is the set of
non-negative integers. The service rate can be changed at the times of arrivals,
departures or phase transitions. Together, the union of these event times and (in a
moment) the added dummy transitions due to uniformization comprise the set of
decision epochs. Based on the queue length, $n$, and the state of the arrival process,
$s$, the controller selects a service rate $\mu_{n,s}$ from the compact set
$\A=[0,\bar{u}]$, $\bar{u}<\infty$. When a service rate $\mu \in \A$ is chosen,
the system incurs a cost at the rate of $c(\mu)$ per unit time.  The cost rate
function, $c(\cdot)$, is defined on $\A$ and is assumed to be strictly convex,
continuously differentiable, strictly increasing and (without loss of generality) such
that $c(0)=0$. Furthermore, a holding or congestion cost is incurred at rate $h(n)$
per unit time when the buffer state is $n$. The holding cost function $h(n)$ is
assumed to be convex, non-decreasing in $n$ and such that $h(0) = 0$ and $\lim_{n\rightarrow\infty}h(n) = \infty$. In the average-cost case we assume $h(\cdot)$ to be a non-decreasing and convex with polynomial rate of growth ($h(n)\leq Cn^p$ for some $C\geq0 , p \in \Z^+$) and again such that $h(0) = 0$ and $\lim_{n\rightarrow\infty}h(n) = \infty$. The assumption about polynomial growth rate of the holding cost function for the average cost case is required for proving the existence of a policy that incurs costs at a finite rate.

Let $\Pi$ be the set of non-anticipating policies. A stationary control policy, $\pi \in
\Pi$, is defined as $\pi=\{\mu(n,s)\mid n\in \N,s\in\mathcal{S}\}$, where
$\mu(n,s)$ is the service rate to be selected when the state of the system is
$(n,s)$. The controller remains idle when the queue is empty i.e, for any policy,
$\mu(0,s)\equiv0$. Thus, given a policy $\pi$, the overall process, $X(t)$, evolves
as a two dimensional continuous time Markov chain on the state space
$\X=\{(n,s)\mid n\in \Z^+,s\in\mathcal{S}\}$. Our objective is to find a control
policy that minimizes the discounted expected cost or average expected cost per
unit time over an infinite time horizon. 

\subsection{The Discounted Expected Cost Formulation}
For $x=(n,s)$ and service rate $\mu \in \A$, let
\begin{align*}
	f(x,\mu):=c(\mu)+h(n).
\end{align*}
Let $\{(X^{\pi}(t),D^{\pi}(t)), t \geq 0 \}$ be the stochastic process representing
the evolution of states and decisions under an admissible policy $\pi$. Given the initial
state $x$ and discount factor $\alpha >0$, the $\alpha-$discounted expected cost
until time $t$ under policy $\pi$ is given by
\begin{align}
	v_{t,\alpha}^\pi(x)&
:=\E_x^{\pi}\left[\int_{0}^{t}e^{-\alpha u}f(X(u),D(u))du\right], \label{eq:total-cost}
\end{align}
where $\E_x[\cdot]:= \E[\cdot|X(0)=x]$. The \textbf{total discounted expected cost} of a policy $\pi$ given that the initial
state of the system is $x$, is
\begin{align*}
	v^\pi_\alpha(x) & :=\lim_{t\to\infty}v^\pi_{t,\alpha}(x).
\end{align*}
The \textbf{optimal total discounted expected cost} is
\begin{align*}
	v_\alpha^\ast(x) & :=\inf_{\pi \in \Pi} v_\alpha^\pi(x).
\end{align*}
A  policy, $\pi^\ast$, is \textbf{total discounted expected cost optimal} if
$v_\alpha^{\pi^\ast}(x)=v_\alpha^\ast(x)$ for all $x \in \X$.

We apply \textit{uniformization} in the spirit of Lippman \cite{Lippman1975} and
consider the discrete time equivalent of the continuous time Markov chain described
above. The uniformization rate is chosen to be $\nu
:=\lambda_L+\bar{\eta}+\bar{u}$, where $\bar{\eta}\geq \max\{-Q_{ss}\mid s\in
\mathcal{S}\}$ is any finite rate larger than the maximum
of the holding time parameters for the phase transition process. 

 Let $v_{\alpha,k}(n,s)$ be the minimum total $\alpha$-discounted expected cost that can be obtained during the last $k$ transitions when starting from state $(n,s)$. Using standard arguments of Markov decision theory \cite{bert}, the discrete-time
finite horizon optimality equations (FHOE) for the system can be written (for each $s
\in \mathcal{S}$):
\begin{subequations}\label{eq:FHOE-simp}
	\begin{align}
		v_{\alpha,k+1}(0,s)&= \frac{1}{\alpha+\nu}\Big[h(0)+\lambda_s v_{\alpha,k}(1,s) \nonumber\\
& \quad +\sum_{s^\prime=1}^{L}
Q_{ss^\prime}v_{\alpha,k}(0,s^\prime)
+(\nu-\lambda_s)v_{\alpha,k}(0,s)\Big] \\
		v_{\alpha,k+1}(n,s)&= \frac{1}{\alpha+\nu}\min_{\mu\in \A}\bigg\{c(\mu)
+h(n)+\mu v_{\alpha,k}(n-1,s)+ \lambda_s v_{\alpha,k}(n+1,s) \nonumber \\
&\quad{}+\sum_{s^\prime=1}^{L}Q_{ss^\prime}v_{\alpha,k}(n,s^\prime)
+\left(\nu-\lambda_s-\mu\right)v_{\alpha,k}(n,s)\bigg\}\quad \text{  for }n \geq 1,
	\end{align}
\end{subequations}
where $v_{\alpha,0}$ is assumed to be zero for each state. Note that the cost
function has compact level sets. That is, $\{((n,s),\mu)| f((n,s), \mu) \leq \beta\}$ is
compact for all $\beta \in \R$. Since the state space is discrete, we may apply
Proposition 3.1 of \cite{Feinberg2007} to get $v_{\alpha,k} \uparrow v_{\alpha}$.
Moreover, $v_{\alpha}$ satisfies \eqref{eq:FHOE-simp} with $v_{\alpha}$ replacing
$v_{\alpha,k}$ on the right hand side and $v_{\alpha,k+1}$ on the left hand side.
The resulting set of equations are called the \textit{discounted cost optimality
equations} (DCOE) and are stated next for later reference (for each
$s\in\mathcal{S}$).
\begin{subequations}\label{eq:DCOE-simp}
	\begin{align}
		v_\alpha(0,s)&= \frac{1}{\alpha+\nu}\left[h(0)+\lambda_s v_\alpha(1,s)+\sum_{s^\prime=1}^{L}
		Q_{ss^\prime}v_\alpha(0,s^\prime)+(\nu-\lambda_s)v_\alpha(0,s)\right] \\
		v_\alpha(n,s)&= \frac{1}{\alpha+\nu}\min_{\mu\in
		\A}\bigg\{c(\mu)+h(n)+\mu v_{\alpha}(n-1,s)+ \lambda_s
		v_\alpha(n+1,s) \nonumber \\
		&\quad{}+\sum_{s^\prime=1}^{L}Q_{ss^\prime}v_\alpha(n,s^\prime)
+\left(\nu-\lambda_s-\mu\right)v_\alpha(n,s)\bigg\}\text{  for }n \geq 1.
	\end{align}
\end{subequations}

\subsection{The Long-Run Average Cost Formulation}
In this section we provide conditions under which an average cost optimal policy
exists and may be computed as a limit of discounted cost optimal policies.
The \textbf{long-run average cost} or \textbf{gain} of a policy $\pi$ given that the
initial state of the system is $x$, is
\begin{align*}
g^{\pi}(x)& :=\limsup_{t\to\infty}v_{t,0}^\pi(x)/t,
\end{align*}
where $v_{t,0}$ is as defined in \eqref{eq:total-cost}. The \textbf{optimal
expected average cost} $g^\ast(x)$ is
\begin{align*}
		g^\ast(x) & :=\inf_{\pi \in \Pi} g^\pi(x),
\end{align*}
and $\pi^\ast$ is an \textbf{average cost optimal policy} if $g^\pi(x)=g^\ast(x)$ for
all $x \in \X$. After uniformization the average cost optimality inequalities (ACOI) (cf.
\cite{sennott89}) are,
\begin{align}\label{eq:ACOI}
	w(n,s)\geq&\frac{1}{\nu}\min_{\mu\in \A}\bigg[-g+c(\mu)+h(n)+
	\lambda_s w(n+1,s)+\mu w((n-1)^+,s) \nonumber\\
	&+\sum_{s^\prime=1}^{L}
	Q_{ss^\prime}w(n,s^\prime)+(\nu-\lambda_s-\mu)w(n,s)\bigg]\text{  for
	}n\geq0, s\in\mathcal{S}.
\end{align}

When the solution, $(w,g)$ to the ACOI exists, $w$ is called a \textit{relative value
function} and $g^\ast(x)=g$ is the optimal long-run expected average cost for any
initial state $x$. 

A solution to the ACOI (\ref{eq:ACOI}) exists under a necessary and sufficient
stability condition which is provided in \eqref{eq:stability} below. This condition
requires that the maximum available service rate is higher than the long-run average
arrival rate and coincides with the one derived by Yechiali \cite{Yechiali} for the
stability of queue with Markov-modulated arrivals. However, since Yechiali used the
balance equations to show the existence of a steady state distribution there is no
guarantee of finite long-run average cost. This is required for the MDP formulation
provided. Since the phase transition process is assumed to be \textit{ergodic}, it has
unique stationary probabilities denoted by, $\{p_1,p_2,\dots,p_L\}$.

\begin{prop}\label{prop:stable}
There exists a stationary policy ,$\pi$, under which the system is stable (steady state
distribution exists) if and only if the maximum available service rate satisfies the
following condition
\begin{equation}
\label{eq:stability}
	\bar{u}>\sum_{s=1}^L p_s\lambda_s.
\end{equation}
Furthermore, the long-run average cost under this policy, $g^{\pi}(x)$, is finite and independent of
the initial state $x$.
\end{prop}
\Pf See Appendix.

In the next proposition, we present results related to the existence of an optimal average-cost policy.  

\begin{prop}\label{prop:sennott}
The following hold
\begin{enumerate}
\item For $\alpha>0$, $v_\alpha(x)$ satisfies the DCOE \eqref{eq:DCOE-simp}.
Moreover, any stationary policy $\pi_\alpha$ that minimizes the right side of
the DCOE \eqref{eq:DCOE-simp} is $\alpha$-discounted expected cost
optimal.
\item If the stability condition (\ref{eq:stability}) holds, we have,
\begin{enumerate}
\item There exists a stationary long-run average expected cost optimal
policy $\pi^\ast=\{\mu^\ast(n,s)\mid n \geq 1,
s\in\{1,2,\dots,L\}\}$  that is a limit of a sequence of
discounted expected cost optimal policies $\{\pi_{\alpha_k}, k
\geq 1\}$. That is,
$\mu^\ast(n,s)=\lim_{k\rightarrow\infty}\mu_{\alpha_k}(n,s)$,
where $\alpha_k\downarrow0$.
\item The long-run average expected cost of policy
$\pi^\ast$ is $g^\ast = \lim_{\alpha\downarrow 0}\alpha
v_\alpha(x)$ for every $x\in\X$. Moreover, there exists a
subsequence $\alpha_k \downarrow 0$ such that $\lim_{k \to \infty}
w_{\alpha_k}(x) :=v_{\alpha_k}(x) - v_{\alpha_k}(\mathbf{0})  =
w(x)$ for a distinguished state $\mathbf{0}$ such that $(w,g^*)$
satisfy the ACOI \eqref{eq:ACOI}.
\end{enumerate}
\end{enumerate}	
\end{prop}
\Pf See Appendix.


\section{Structural Properties of Optimal Policies}\label{sec:structure}
In this section we derive structural results for optimal policies for both the discounted
cost and the average cost criterion.
In a manner similar to \cite{George}, we use following definitions to simplify the
optimality equations,
\begin{align}
	\label{eq:CONVCONJ}
	y_{\alpha}(0,s)&=0\text{ for } s=1,2,\dots,L,\\\nonumber
	y_{\alpha}(n,s)&=v_{\alpha}(n,s)-v_{\alpha}(n-1,s) \text{ for } n\in\N,s\in\mathcal{S},\\\nonumber
	\phi(y)&=\max_{\mu \in \A}\{\mu y-c(\mu)\}, \text{ and } \\\nonumber
	\psi(y)& = \arg\max_{\mu \in \A}\{\mu y-c(\mu)\},
\end{align}
where the argmax is a singleton by the assumptions on $c$ (cf. Section 4.3 of
\cite{Ata2005}). The definitions above yield the following simplified form of the
DCOE \eqref{eq:DCOE-simp}:
\begin{align}\label{eq:DCOE}
	v_\alpha(n,s)&=\frac{1}{\alpha+\nu}\bigg[h(n)-\phi(y_{\alpha}(n,s))+
	\lambda_s v_\alpha(n+1,s)+\sum_{s^\prime=1}^{L}
	Q_{ss^\prime}v_\alpha(n,s^\prime)  \nonumber\\
	&\quad{}+(\nu-\lambda_s)v_\alpha(n,s)\bigg]\text{  for
	}n\in\Z^+,s\in\mathcal{S}.
\end{align}

In order to derive structural results for an optimal discounted expected cost policy,
we make use of several important properties of functions $\phi(y)=\max_{x \in
\A}\{yx-c(x)\}$ and its associated maximizers $\psi(y) = \arg\max_{x \in
\A}\{yx-c(x)\}$ that were introduced in the DCOE \eqref{eq:DCOE-simp}. Recall
that the \textit{conjugate} of $c(\cdot)$, $\phi(\cdot)$, is convex (cf. \cite{Boyd}).
Moreover, $\psi(y)$ is continuous, non-decreasing and equals $\phi^\prime(y)$
wherever the derivative exists. As described in \cite{Ata2005}, since
$(c^\prime)^{-1}(\cdot)$ is well-defined, continuous and strictly increasing we
have the following characterization of the function $\psi(\cdot)$
\begin{align}
	\psi(y) & =
		\begin{cases}
		  0  & \text{ if } y \leq c^{\prime}(0),\\
		      (c^\prime)^{-1}(y) &  \text{ if } c^{\prime}(0)< y <
		      c^{\prime}(\bar{u}), \\
		  \bar{u} & \text{ if }  y >c^{\prime}(\bar{u}).
		\end{cases}
\end{align}
It may also be established (see \cite{George}) that $\phi(\cdot)$ is continuous and
non-decreasing with the following characterization
\begin{align}\label{eq:phi}
\phi(y) & =\begin{cases}
		0  & \text{ if } y < 0,\\
		\int_0^y\psi(x)dx   & \text{ if } y \geq 0.\\
		\end{cases}
\end{align}
\subsection{Monotone In The Number of Customers}
We show the intuitive result that there exists an optimal policy that is monotone in $n$. We
note that the structural part of the result could also be proven via the event-based
dynamic programming framework of Koole \cite{Koole1998}. We provide what we believe
is an equally simple proof here for completeness.
\begin{prop}\label{prop:mono-values}
The following hold
\begin{enumerate}
\item For each $s\in\mathcal{S}$, the optimal discounted expected cost
value function, $v_\alpha(n,s)$, satisfies the DCOE \eqref{eq:DCOE-simp}
and is a non-decreasing, convex function of $n$.
\item There exists a discounted expected cost optimal policy
$\{\mu_{\alpha}(n,s), n \geq 1, s \in \mathcal{S}\}$ that is
non-decreasing in $n$ for each $s\in\mathcal{S}$.
\item Under the assumptions that the holding cost is non-decreasing and convex with polynomial rate of growth and
\eqref{eq:stability} hold, there exists a long-run average optimal policy,
$\{\mu(n,s), n \geq 1, s \in \mathcal{S}\}$ that is non-decreasing in $n$
for each $s\in\mathcal{S}$.
\end{enumerate}
\end{prop}
\Pf 
%
We use induction and the FHOE \eqref{eq:FHOE-simp} to prove the first result. The
result holds trivially for $k=0$. For the inductive step, suppose $v_{\alpha,k}(\cdot
,s)$ is non-decreasing and convex on $\Z^+$ for each $s\in\mathcal{S}$. Let
$u_n=\mu_{\alpha,k}(n,s)$ be the optimal service rate for the $(k+1)$-stage
problem when the state is $(n,s)$. Suppose we use the potentially sub-optimal
decision $u_n$ when the state is $(n-1,s)$. The FHOE \eqref{eq:FHOE-simp} yield,
\begin{align*}
	v_{\alpha,k+1}(n-1,s)&\leq \frac{1}{\alpha+\nu}\bigg[c(u_n)+h(n-1)
+u_n v_{\alpha,k}((n-2)^+,s)+ \lambda_s v_{\alpha,k}(n,s)\\\nonumber
	&\quad+\sum_{s^\prime=1}^{L}Q_{ss^\prime}v_{\alpha,k}(n-1,s^\prime)
	+\left(\nu-\lambda_s-u_n\right)v_{\alpha,k}(n-1,s)\bigg],\\
	&\leq \frac{1}{\alpha+\nu}\bigg[c(u_n)+h(n)+u_n v_{\alpha,k}(n-1,s)+ \lambda_s
	v_{\alpha,k}(n+1,s)\\\nonumber
	&\quad+\sum_{s^\prime \neq s} Q_{ss^\prime}v_{\alpha,k}(n,s^\prime)
	+\left(\nu + Q_{ss} - \lambda_s-u_n\right)v_{\alpha,k}(n,s)\bigg]\\
	&=v_{\alpha,k+1}(n,s),
\end{align*}
where the second inequality follows from the induction hypothesis. Thus,
$v_{\alpha,k}$ is non-decreasing for all $k$.

To show convexity note that by the inductive hypothesis,
$y_{\alpha,k}(n+1,s)=v_{\alpha,k}(n+1,s)-v_{\alpha,k}(n,s)$ is a non-decreasing
function of $n$ for each $s\in\mathcal{S}$. Let $u_{n+1}=\mu_{\alpha,k}(n+1,s)$
be the optimal rate for the $(k+1)$-stage problem when the state is $(n+1,s)$ and
$u_{n-1}=\mu_{\alpha,k}(n-1,s)$ be the optimal rate when the state is $(n-1,s)$.
The DCOE \eqref{eq:DCOE-simp} imply (for $n\geq1$)
\begin{align*}
	(\alpha+\nu)y_{\alpha,k+1}(n+1,s)&\geq
		  h(n+1)-h(n)-u_{n+1}(y_{\alpha,k}(n+1,s)-y_{\alpha,k}(n,s))\\\nonumber
	   &\quad +\lambda_s y_{\alpha,k}(n+2,s)+\sum_{s^\prime=1}^{L}Q_{s,s^\prime}y_{\alpha,k}(n+1,s^\prime).\\\nonumber
\end{align*}
Similarly for $y_{\alpha,k}(n,s)=v_{\alpha, k}(n,s)-v_{\alpha,k}(n-1,s)$,
\begin{align*}
	(\alpha+\nu)y_{\alpha, k+1}(n,s)&\leq
		h(n)-h(n-1)+(\nu-\lambda_s)y_{\alpha,k}(n,s)\\
	&\quad -u_{n-1}(y_{\alpha,k}(n,s)-y_{\alpha,k}(n-1,s))+\lambda_s y_{\alpha,k}(n+1,s)\\
	&\quad +\sum_{s^\prime=1}^{L}Q_{s,s^\prime}y_{\alpha,k}(n,s^\prime).
\end{align*}
Using the definitions of $\nu=\lambda_L+\bar{\eta}+\bar{u}$ and
$\bar{\textbf{Q}}=\bar{\eta}\textbf{I}+\textbf{Q}$ we have,
\begin{align*}
	(\alpha+\nu)(y_{\alpha,k+1}(n+1,s)-y_{\alpha,k+1}(n,s))&\geq
h(n+1)-2h(n)+h(n-1)\\
	   &\quad +(\lambda_L+\bar{u} -\lambda_s-u_{n+1})\left(y_{\alpha,k}(n+1,s)-y_{\alpha,k}(n,s)\right)\\
	   &\quad +u_{n-1}(y_{\alpha,k}(n,s)-y_{\alpha,k}(n-1,s))\\
	   &\quad +\lambda_s\left(y_{\alpha,k}(n+2,s)
	       -y_{\alpha,k}(n+1,s)\right)\\
&\quad + \sum_{s^\prime =1}^L \bar{Q}_{s,s^\prime}
\left(y_{\alpha,k}(n+1,s^\prime)-y_{\alpha,k}(n,s^\prime)\right)\\
&\geq 0.
\end{align*}
The second inequality follows as $h(n)$ is convex, the coefficients of
$y_{\alpha,k}$ terms are non-negative and the inductive hypothesis. So
$y_{\alpha,k}(\cdot,s)$ is non-decreasing on $\Z^+$ for all $s\in\mathcal{S}$ as
required. Taking limits as $k \to \infty$ yields that $v_{\alpha}$ is non-decreasing
and convex; the first result is proven.

Since the function $\psi(\cdot)$ is non-decreasing and $\mu_{\alpha}(n,s)=\psi(y_{\alpha}(n,s))$, we conclude that there exists an optimal policy for the discounted cost problem that is
monotonically nondecreasing in the queue length for each $s \in \mathcal{S}$. This is the second result.

For the third result, consider a subsequence of discount factors $\{\alpha_i, i \geq 0\}$ such that $\alpha_i \to 0$ and corresponding
discounted cost optimal policies $\mu_{\alpha_i}(\cdot,\cdot)$ that
converge to an average cost optimal policy $\mu(\cdot, \cdot)$ (see Proposition \ref{prop:sennott}). The previous result
implies that for each fixed $s \in \mathcal{S}$, $\mu_{\alpha_i}(n,s) \leq
\mu_{\alpha_i}(n+1,s)$. Thus, the same inequality holds for $\mu(\cdot,\cdot)$. \EndPf

We remark that we have explicitly used the fact that argmax in $\psi$ is a singleton
(which follows from the strict convexity assumption on $c(\cdot)$). When the
convexity is not assumed to be strict, the results still hold, but we need to take
care to define $\psi$ as the minimal element of the argmax and consider a
subsequence of discount factors such that $w_{\alpha_i}(x) = v_{\alpha_i}(x) -
v_{\alpha_i}(\mathbf{0}) \to w(x)$. Since $v_{\alpha_i}(n,s)$ is non-decreasing
in $n$, so is $w$ and proof in the average cost case follows in the same way as the
discounted cost case except that we use the ACOI instead of the DCOE.

\subsection{Monotone in the phase process}
Since the states of the phase transition process are ordered such that
$\lambda_1\leq\lambda_2\leq\dots\leq\lambda_L$, one might conjecture that the
optimal policy is non-decreasing in the phase state, $s$, for each congestion level,
$n$. However, we present two examples to show that depending on the transition
structure of the phase process, this property may \textit{not} hold. In both
examples, we use value iteration with $\alpha=0.05$ to compute the optimal policy
numerically. We consider an exponential cost rate function, $c(\mu)=e^{\mu}-1$,
and a linear holding cost function, $h(n)=n$. The set of permissible service rates is $\A=[0,5]$. 

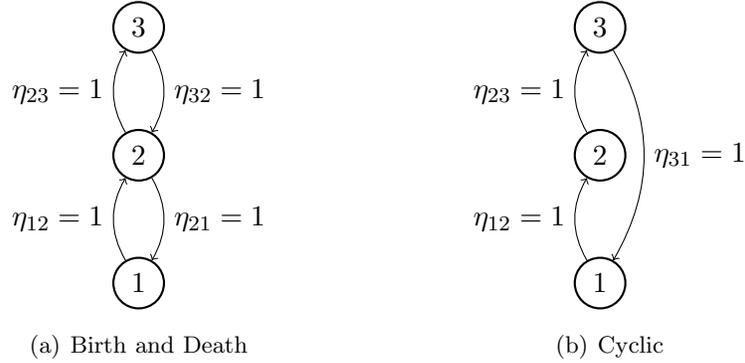
\begin{figure}[htbp]
\centering
\subfigure[Birth and Death]{
		\begin{tikzpicture}[auto]
		\node (top) [circle,thick,draw=black,minimum size=5mm]{$3$};
		\node (mid) [circle,thick,draw=black,minimum size=5mm][below=of top]
		{$2$};
		\node (low) [circle,thick,draw=black,minimum size=5mm][below=of
		mid]{$1$};
		\draw [->] (low) to [bend left=30] node{$\eta_{12}=1$}(mid);
		\draw [->] (mid) to [bend left=30]node{$\eta_{21}=1$}(low);
		\draw [->] (mid) to [bend left=30] node{$\eta_{23}=1$}(top);
		\draw [->] (top) to [bend left=30]node{$\eta_{32}=1$}(mid);
		\end{tikzpicture}
		\label{fig:rw}
}       \hspace{20mm}
\subfigure[Cyclic]{
		\begin{tikzpicture}[auto]
		\node (top) [circle,thick,draw=black,minimum size=5mm]{$3$};
		\node (mid) [circle,thick,draw=black,minimum size=5mm][below=of top]
		{$2$};
		\node (low) [circle,thick,draw=black,minimum size=5mm][below=of
		mid]{$1$};
		\draw [->] (low) to [bend left=30] node{$\eta_{12}=1$}(mid);
		\draw [->] (mid) to [bend left=30] node{$\eta_{23}=1$}(top);
		\draw [->] (top) to [bend left=30]node{$\eta_{31}=1$}(low);
		\end{tikzpicture}
		\label{fig:cyc}
}
\caption{Transition structure of Phase process for Examples \ref{ex:mono} and \ref{ex:non-mono}.}
\end{figure}

\begin{ex}
\label{ex:mono}
In this example, the phase process is a birth and death process on the states
$\{1,2,3\}$. See Figure \ref{fig:rw}. The infinitesimal generator for the phase
process is given by
\[\mathbf{Q}=\left[\begin{matrix}-1&1&0\\1&-2&1\\0&1&-1\end{matrix}\right],\]
and arrival rates are $\lambda_1=0.5,\lambda_2=1$ and $\lambda_3=1.25$.

\begin{figure}
\centering
	\includegraphics[totalheight=0.4\textheight]{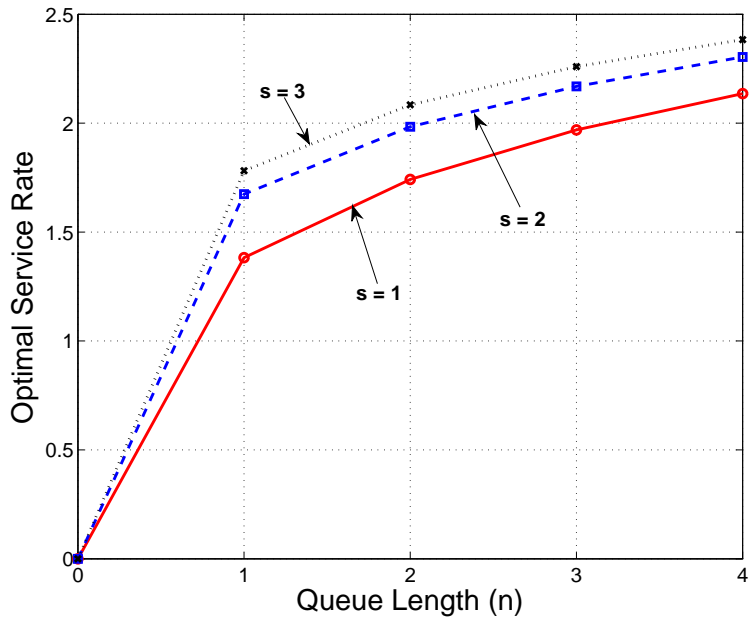}
	\caption{Structure of Optimal Policy for Example 3.1}
	\label{fig:ex31}
\end{figure}
\begin{figure}
\centering
    \includegraphics[totalheight=0.4\textheight]{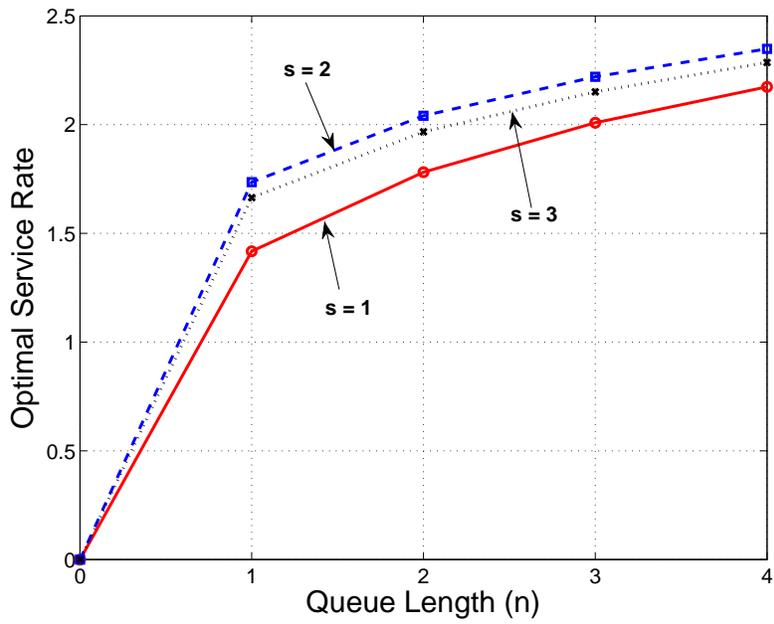}\label{fig:non-mono1}
    \caption{Structure of Optimal Policy for Example 3.2}
    \label{fig:ex32}
\end{figure}

Figure \ref{fig:ex31} shows that the optimal policy is a non-decreasing function of $n$ for each $s$.
It should also be clear that the optimal service rates are non-decreasing in $s$ for each $n$.

\end{ex}

\begin{ex}\label{ex:non-mono}

Consider a phase process with a cyclic transition structure on the set of states
$\{1,2,3\}$. See Figure \ref{fig:cyc}. The infinitesimal generator matrix for the
phase process is
\[\mathbf{Q}=\left[\begin{matrix}-1&1&0\\0&-1&1\\1&0&-1\end{matrix}\right]\]
and the arrival rates are $\lambda_1=0.5,\lambda_2=1$ and $\lambda_3=1.25$.
Figure \ref{fig:non-mono1} shows that the optimal policy is non-decreasing in $n$ for
each $s$. However, it is clear from this figure that the service rates
are not monotone in $s$ when queue length is 4.

\end{ex}

In Example \ref{ex:non-mono}, the phase process transitions from the highest arrival
intensity state, 3, to the lowest arrival intensity state, 1. This causes the optimal service rate to be higher in state 2 as compared to state 3 for some congestion levels and thereby renders an optimal policy that is not
monotone in $s$. These examples beg the question, is there a reasonable assumption
under which the optimal policy is monotone in $s$? \textit{Stochastic monotonicity}
of the phase transition process is one such assumption.

\subsubsection{Stochastic Monotonicity for Continuous Time Markov Chains}
Intuitively, stochastic monotonicity means that given the arrival process is in a high
arrival intensity state, the future states it will encounter are in some sense worse
(in terms of arrival intensity) than if the process is in a low arrival intensity state.
This leads to the following definitions (see for example  Keilson and Kester
\cite{Keilson1977}).
\begin{defn}
Given two probability vectors \textbf{p} and \textbf{q}, a stochastic matrix
\textbf{M} and a homogeneous Markov chain $\{X(t),t\geq0\}$ with probability
transition function $\textbf{P}(t)=\textbf{P}_{ij}(t)$.
\begin{enumerate}
\item \textbf{p} stochastically dominates \textbf{q}
(\textbf{p}$\geq_{st}$\textbf{q}) iff
$\sum_{i=n}^Np_i\geq\sum_{i=n}^N q_i$, $n=1,2,\dots$.
\item Letting $\mathbf{M}^i$ denote the $i^{th}$ row of the matrix,
\textbf{M} is called stochastically monotone if
$\mathbf{M}^k\geq_{st}\mathbf{M}^l$, whenever $k>l$.
\item $\{X(t),t\geq0\}$ is said to be \textit{stochastically monotone} if
$\textbf{P}(t)$ is monotone.
\end{enumerate}
\end{defn}
Note that the transition structure shown in Example \ref{ex:mono} is stochastically
monotone while that for Example \ref{ex:non-mono} is not (this is trivial to see once the underlying Markov Chain is uniformized). Some useful stochastic processes that have the stochastic monotonicity property include the birth-death process, the simple random walk, the age of renewal process with decreasing failure rate \cite{Keilson1977}. In particular, the simple 2-state MMPP fluctuating between high arrival rate and low arrival rate considered by Gupta \cite{Gupta2006} and Shah-Heydari \cite{Shah-Heydari1998}, is also stochastically monotone. The following provides alternative methods for specifying when a transition
matrix is stochastically monotone (again refer to Keilson and Kester
\cite{Keilson1977}).
\begin{prop}\label{prop:tfae}
For monotone matrices the following are equivalent.
\begin{enumerate}
\item \textbf{M} is monotone.
\item $(\mathbf{T}^{-1}\mathbf{M}\mathbf{T})_{ij}\geq0$ where
\textbf{T} is a square matrix with 1's on or below the diagonal.
\item \textbf{pM}$\geq_{st}$\textbf{qM} for all probability vectors
\textbf{p},\textbf{q} with \textbf{p}$\geq_{st}$\textbf{q}.
\item \textbf{Mv} is non-decreasing for all non-decreasing vectors
\textbf{v}.
\end{enumerate}
\end{prop}
Since for a continuous-time Markov chain with generator matrix $\mathbf{Q}$, stochastic
monotonicity implies that the generator for the phase process satisfies the property
$(\textbf{T}^{-1}\textbf{Q}\textbf{T})_{ij}\geq0, i\neq j$ where
\textbf{T} is a square matrix with 1's on or below the diagonal(\cite{Keilson1977}). Choosing a
uniformizing constant $\bar{\eta}$ yields
$(\textbf{T}^{-1}(\bar{\eta}\textbf{I}+\textbf{Q})\textbf{T})\geq0$ in all
elements. Thus using the second and fourth parts of Proposition \ref{prop:tfae} we have that \textbf{Q} is such
that $\bar{\textbf{Q}}:=\bar{\eta}\textbf{I}+\textbf{Q}$ satisfies the property
that $\bar{\textbf{Q}}\textbf{v}$ is non-decreasing for all non-decreasing vectors
\textbf{v}. This leads to the next result; the main result of this section.

\begin{thm}
Suppose that the phase transition process is stochastically monotone.
\begin{enumerate}
\item For each $n\in \Z^+$, $y_{\alpha,k}(n,s)$ is non-decreasing function
of $s$.
\item There exists a discounted cost optimal policy, $\mu_{\alpha}(n,s)$,
that is non-decreasing in $s$ for each $n$.
\item Under the assumptions that the holding cost is non-decreasing and convex with polynomial rate of growth and \eqref{eq:stability} hold, there exists an average cost optimal policy $\mu(n,s)$, that is non-decreasing in $s$ for each $n$.
\end{enumerate}
\end{thm}
\Pf We show the first result by induction. The second and third results follow in an
analogous manner to Proposition \ref{prop:mono-values}. The statement holds
trivially for $k=0$. Assume it holds for $k$. Using the definitions
$\nu=\lambda_L+\bar{\eta}+\bar{u}$ and
$\bar{\textbf{Q}}=\bar{\eta}\textbf{I}+\textbf{Q}$ we have for $s>1$
\begin{align}\label{eq:3.4}
	\lefteqn{(\alpha+\nu)\left(y_{\alpha,k+1}(n+1,s)-y_{\alpha,k+1}(n+1,s-1)\right)} \nonumber \\
&\quad =\phi\left(y_{\alpha,k}(n,s)\right)-\phi\left(y_{\alpha,k}(n,s-1)\right)
+\lambda_{s}\left(y_{\alpha,k}(n+2,s)-y_{\alpha,k}(n+2,s-1)\right) \nonumber\\
& \qquad +(\lambda_s-\lambda_{s-1})\left(y_{\alpha,k}(n+2,s-1)
-y_{\alpha,k}(n+1,s-1)\right) \nonumber \\
	&\qquad +(\lambda_L-\lambda_s)\left(y_{\alpha,k}(n+1,s)
-y_{\alpha,k}(n+1,s-1)\right) \nonumber \\
	&\qquad+\left[\bar{u}\left(y_{\alpha,k}(n+1,s)-y_{\alpha,k}(n+1,s-1)\right)\right. \nonumber\\
	&\qquad \left.-\left(\phi(y_{\alpha,k}(n+1,s))-\phi\left(y_{\alpha,k}(n+1,s-1)\right)\right)\right]\nonumber \\
	&\qquad +\bigg[\sum_{s^\prime=1}^{L}\bar{Q}_{s,s^\prime}y_{\alpha,k}(n+1,s^\prime)
-\sum_{s^\prime=1}^{L}\bar{Q}_{s-1,s^\prime}y_{\alpha,k}(n+1,s^\prime)\bigg]
\end{align}
Note that the inductive hypothesis implies that the first four terms in the RHS of \eqref{eq:3.4} are non-negative. Now as $\phi(y)=\int_0^{y}\psi(x)dx$ and $\psi(y)\leq\bar{u}$,
\begin{align*}
	 \phi(y_{\alpha,k}(n+1,s))-\phi(y_{\alpha,k}(n+1,s-1))&=\int_{y_{\alpha,k}(n+1,s-1)}^{y_{\alpha,k}(n+1,s)}\psi(x)dx\\
	&\leq \bar{u}(y_{\alpha,k}(n+1,s)-y_{\alpha,k}(n+1,s-1)).
\end{align*}
So the next to last term in RHS of \eqref{eq:3.4} is non-negative. Furthermore, since the inductive hypothesis and the assumption on $\bar{Q}$ implies the last term in the RHS of \eqref{eq:3.4} is non-negative. Thus, we conclude from the
induction hypothesis that
\begin{equation*}
	(\alpha+\nu)(y_{\alpha,k+1}(n+1,s)-y_{\alpha,k+1}(n+1,s-1)) \geq 0,
\end{equation*}
as desired.
\EndPf
%
%


\section{Numerical Study}\label{sec:numeric}
This section provides two insights using numerical examples. First, we compare the optimal control policy with
two natural heuristics. When the environment is changing, it seems a
decision-maker that is not armed with the current research might take one of two courses.
(S)he might choose to ignore the state change of the environment
altogether, or she might treat each state change as permanent and react
accordingly. In either case, the resulting control policies are heuristics when
compared to the optimal control that takes into account both the phase and queue
length processes. The first goal is then to compare the optimal policy with these
heuristics. Second, as alluded to in Section \ref{sec:intro}, the current model can
act as a heuristic itself when compared to a model with NHPP arrivals. We analyze
when this is a reasonable approximation.

\subsection{Comparison with Heuristics}
 In this section we present a comparison of the performance of optimal policy with several heuristics. Furthermore, we compare the optimal cost achievable with state-dependent service rates with the optimal cost achievable when the service rate is fixed for all states. As mentioned in George and
Harrison \cite{George}, the difference between these costs represent the economic value of a responsive mechanism.  The first heuristic that we consider uses the optimal control for an average cost
problem where the arrival process is Poisson with the long-run mean arrival rate of
the MMPP. When applied to the original model, this policy is a function of the queue
length only. We call this heuristic the \textit{Average Rate Method} (ARM).

Since the state of the arrival process is known, the decision-maker may solve the
stationary model with each potential arrival rate and change the service rate
according to the current state of the arrival process. That is to say, a second
heuristic is derived in the following way:
\begin{enumerate}
\item Compute the service rate control average cost optimal policy, $\pi^h_s$,
for a system with Poisson arrivals with rate $\lambda_s$ for each intensity
level $s\in\mathcal{S}$.
\item The heuristic for the Markov-Modulated queue is obtained by using
$\pi^h_s$ when the state of the process is $(n,s)$.
\end{enumerate}
This heuristic is referred to as the \textit{Phase Rate based Method} (PRM). Note
that the long-run average arrival rate used in computing the ARM policy is
influenced by both the infinitesimal generator matrix and the arrival rates of the
phase process while the PRM policy relies only on the arrival rates.

When the service rate is fixed for all states (open loop policy), the queue operates as an MMPP/M/1 queue. For a given service rate, the transition matrix and corresponding steady state distribution can be computed numerically. Based on the steady state distribution one can determine the long-run average cost corresponding to that service rate. We then use a 1-D search procedure to find the service rate that minimizes long-run average cost.  This is called the \textit{Fixed Rate} policy.

A numerical study comparing the performance of these heuristics for various test
cases is provided in Examples \ref{ex:numerics1} and \ref{ex:numerics2}. In all
cases, the policies and average cost are computed using value iteration where the
queue length is truncated at 50. We use the cost rate function
$c(\mu)=e^{\mu}-1$ and holding cost rate $h(n)=n$. Service rates are allowed to
be chosen from $\A=[0,15]$. In each case the arrival rates change in accordance
with the phase state $\{1, 2, \ldots, 8\}$ with the arrival rates as shown in Table
\ref{tab:par_ex1}.

\begin{table}[ht!]
\centering
\begin{tabular}{|c|c|c|c|c|c|c|c|c|}
\hline
	  	& \multicolumn{8}{c|}{ Arrival Rate in Phase State}\\
\hline
Case      & 1     & 2     & 3     & 4     & 5     & 6     & 7     & 8 \\\hline
I & 0.1   & 0.35  & 0.6   & 0.85  & 1.1   & 1.35  & 1.6   & 1.85 \\
II & 0.1   & 0.6   & 1.1   & 1.6   & 2.1   & 2.6   & 3.1   & 3.6 \\
III & 0.1   & 0.85   & 1.6   & 2.35   & 3.1   & 3.85   & 4.6   & 5.35 \\
\hline
\end{tabular}%
\caption{Arrival Rate Parameters for Phase Transition Process in Examples
\ref{ex:numerics1} and \ref{ex:numerics2}}
\label{tab:par_ex1}%
\end{table}%

\begin{ex}
\label{ex:numerics1}

\noindent Suppose that phase process is a birth and death process on states
$\{1,2,\dots,8\}$ (recall Figure \ref{fig:rw}). Fix $c>0$. The transition rates for the
phase process are $\eta_{i,i+1}=\eta_{i,i-1}=c$ for $2\leq i\leq 7$,
$\eta_{1,2}=c$ and $\eta_{8,7}=c$. A higher value for $c$ means that the
phase process transitions faster between the arrival phases. We refer to $c$ as
the \textit{fluctuation rate scaling parameter}. For numerical analysis,
we consider three sets of arrival rates for the phase process as shown in Table
\ref{tab:par_ex1}. For each set, the parameter $c$ takes values $0.25,0.50,0.75$
and $1.00$ resulting in a total of 12 different scenarios for the arrival process. Table \ref{tab:ex_41} shows the results for all heuristics.

\begin{table}[hbt!]
\centering
\begin{tabular}{|cc|c|c|c|c|}
\hline
\multicolumn{2}{|c|}{Scenarios} & \multicolumn{4}{c|}{Gain (\% Sub-Optimal)}  \\
\hline
Arrival Rates & $c$     & Optimal & ARM   & PRM & Fixed Rate \\\hline
Case I & 0.25  & 4.3651 & 4.4650 (2.29 \%) & 4.3676 (0.06 \%)  & 7.6841 (76.03 \%)\\
& 0.50   & 4.3196 & 4.3974 (1.80 \%) & 4.3254 (0.13 \%) & 7.3185 (69.43 \%)\\
& 0.75  & 4.2818 & 4.3455 (1.49 \%) & 4.2909 (0.21 \%)& 7.0223 (64.01 \%)\\
& 1.00     & 4.2494 & 4.3031 (1.27 \%) & 4.2618 (0.29 \%)& 6.8399(60.96 \%) \\\hline
Case II & 0.25  & 15.5713 & 16.9349 (8.76 \%) & 15.7936 (1.43 \%)& 24.8200 (59.41 \%) \\
& 0.50   & 14.8674 & 15.6939 (5.56 \%) & 15.2599 (2.64 \%)& 22.5509 (51.68\%) \\
& 0.75  & 14.3638 & 14.9444 (4.04 \%) & 14.8821 (3.61 \%)& 21.1000(46.9\%) \\
& 1.00     & 13.9776 & 14.4189 (3.16 \%) & 14.5924 (4.40 \%)& 20.1360(44.06\%) \\\hline
Case III & 0.25  & 47.6797 & 51.9918 (9.04 \%) & 49.6854 (4.21 \%)& 61.1588(28.27\%) \\
& 0.50   & 42.3561 &44.4741 (5.00 \%) & 45.7978 (8.13 \%)& 55.8678 (31.9\%) \\
& 0.75  & 39.2816 & 40.6579 (3.51 \%) & 43.7541 (11.39 \%)& 51.8600 (32.02\%) \\
& 1 .00    & 37.2150 & 38.2310 (2.73 \%) & 42.3809 (13.88 \%)& 48.9160 (31.44\%) \\
\hline
\end{tabular}%
\caption{Average Cost Rates and Percentage Difference between Optimal and
Heuristic Policies for Example \ref{ex:numerics1}.}
\label{tab:ex_41}%
\end{table}%

A few observations are in order. Ignoring the dynamic state information of the
phase process (and using the ARM policy) is more costly when the fluctuation
parameter is lower. This stands to reason since the phase process can be in a state
for a long period of time, while the ARM policy assumes the arrival rate is the mean
arrival rate. In Case III, when the arrival rate change is the most between phase
states and with the slowest rate of changing states, the percent sub-optimality for ARM is
above 9\%. If we try to approximate the state changes with stationary processes
(using PRM) we see that again, the percent sub-optimality is high (above 13\%) but
this time when the fluctuation parameter is highest.

\begin{figure}[ht!]
\centering
\subfigure[Small Difference in Phase Intensities (Case I)]
{\includegraphics[totalheight=0.28\textheight]{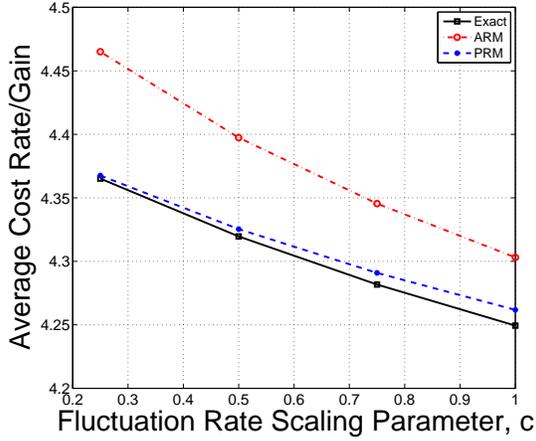} \label{fig:small-diff}}
\subfigure[Med. Difference in Phase Intensities (Case II)]
{\includegraphics[totalheight=0.28\textheight]{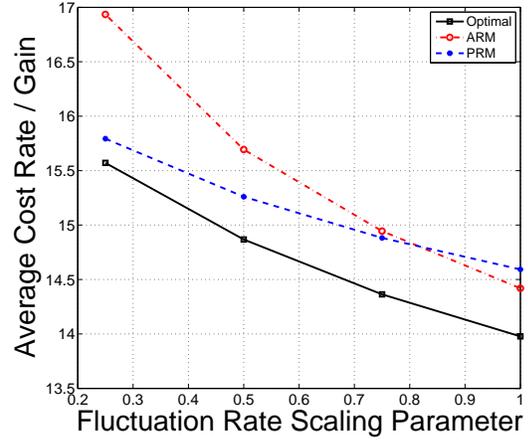} \label{fig:med-diff}}
\subfigure[Large Difference in Phase Intensities (Case III)]
{\includegraphics[totalheight=0.28\textheight]{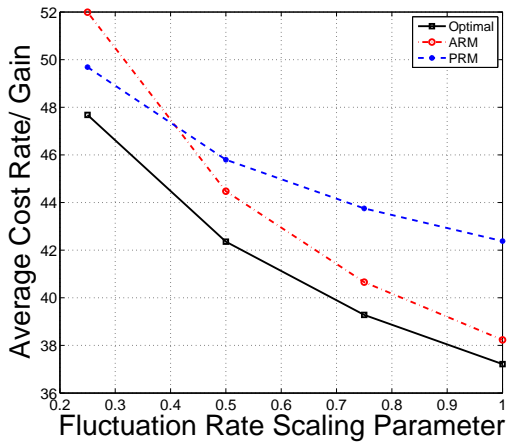} \label{fig:large-diff}}
\subfigure[Heuristic Vs Opt. Policy (Case III, $n = 2$)]
{\includegraphics[totalheight=0.27\textheight]{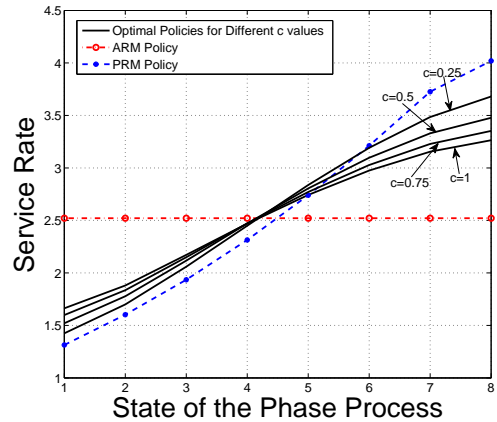} \label{fig:heur-v-opt}}
\caption{Comparison of Gain Values for Heuristic and Optimal Policies for Example \ref{ex:numerics1}}.
\label{fig:ex41}
\end{figure}
Figures \ref{fig:small-diff}-\ref{fig:large-diff} show the change in the average cost
rate under the heuristic and optimal policies as a function of the parameter $c$ for
the three arrival cases. Figure \ref{fig:heur-v-opt} shows a comparison of the
heuristics and optimal policies for the arrival rates in Case III and $n=2$ for various
values of fluctuation rate parameter (the behavior is similar for other values of
$n$). It should be clear that the PRM policy outperforms the ARM policy in most
cases except for high values of $c$ in Case III. Moreover, one should note that the
performance of the ARM policy improves while that of PRM policy degrades in
comparison to the optimal policy as the fluctuation rate parameter increases.
Intuitively this seems reasonable since  for low values of $c$ the phase process
spends more time in each phase. Therefore the PRM policy, that in each phase
applies the optimal policy for a stationary M/M/1 queue with arrival rate of that
particular phase, performs better than the ARM policy. In fact, if $c=0$, the PRM
policy is optimal since the phase process is stationary with the arrival rate of initial
phase.

At high values of $c$ (since the system sees more and more transitions), the arrival
process behaves like a Poisson process with the average arrival rate of MMPP.
Therefore the PRM policy gets penalized more in comparison to the ARM policy in
this case. In Figure \ref{fig:heur-v-opt} we also see that the change in the optimal
service rates as a function of phase state for each congestion level is lower for
higher values of $c$.

\end{ex}


\begin{ex}
\label{ex:numerics2}

\noindent In this example we study a cyclic phase process (cf. Figure \ref{fig:cyc})
on the states $\{1,2,\dots,8\}$.  The transition rates for the phase process are
$\eta_{i,i+1}=c$ for $1\leq i\leq 7$ and $\eta_{8,1}=c$. Similar to the previous
example we perform the numerical analysis for 12 scenarios for the phase process;
three different sets of arrival rates given in Table \ref{tab:par_ex1} and for each
set of arrival rates, $c$ takes values $0.25,0.50,0.75$ and $1.00$.
\begin{table}[h!]
\centering
\begin{tabular}{|cc|c|c|c|c|}
\hline
\multicolumn{2}{|c|}{Scenarios} & \multicolumn{4}{c|}{Gain (\% Sub-Optimal)}  \\
\hline
Arrival Rates & $c$     & Optimal & ARM   & PRM & Fixed Rate\\\hline
Case I & 0.25  & 4.1872 & 4.2295 (1.01 \%) & 4.2267 (0.94 \%)& 6.3440 (51.51\%) \\
& 0.50   & 4.0603 & 4.085 (0.61 \%) & 4.1204 (1.48 \%) & 5.9620 (46.84\%)\\
& 0.75  & 3.988 & 4.0051 (0.43 \%) & 4.0574 (1.73 \%)& 5.7700 (44.68\%) \\
& 1.00     & 3.9423 & 3.9549 (0.32 \%) & 4.0166 (1.89 \%)& 5.6647 (43.69\%) \\\hline
Case II & 0.25  & 12.894 & 13.2042 (2.41 \%) & 13.9767 (8.39 \%)& 17.2439 (33.74 \%) \\
& 0.50   & 11.9656 & 12.1319 (1.39 \%) & 13.2268 (10.54 \%)& 15.6149 (30.5 \%) \\
& 0.75  & 11.5435 & 11.6531 (0.95 \%) & 12.8573 (11.38 \%)& 14.9350 (29.38 \%) \\
& 1.00     & 11.2996 & 11.3786 (0.70 \%) & 12.6374 (11.84 \%)& 14.51 (28.43 \%) \\\hline
Case III & 0.25  & 31.2724 & 32.1887 (2.93 \%) & 39.4752 (26.23 \%)& 35.5711 (13.75 \%) \\
& 0.50   & 28.3046 & 28.7893(1.71 \%) & 37.1449 (31.23 \%)& 33.8800 (19.7 \%) \\
& 0.75  & 27.0506 & 27.3664 (1.16 \%) & 36.0660 (33.33 \%)& 32.7185 (20.95\%) \\
& 1.00    & 26.3445 & 26.5702 (0.86 \%) & 35.4401 (34.53 \%)& 31.9436 (21.25 \%) \\
\hline
\end{tabular}%
\caption{Average Cost Rates and Percentage Difference between Optimal and
Heuristic Policies for Example \ref{ex:numerics2}}
\label{tab:ex_42}%
\end{table}%

\begin{figure}[ht!]
\centering
\subfigure[Small Difference in Phase Intensities (Case I)]
{\includegraphics[totalheight=0.28\textheight]{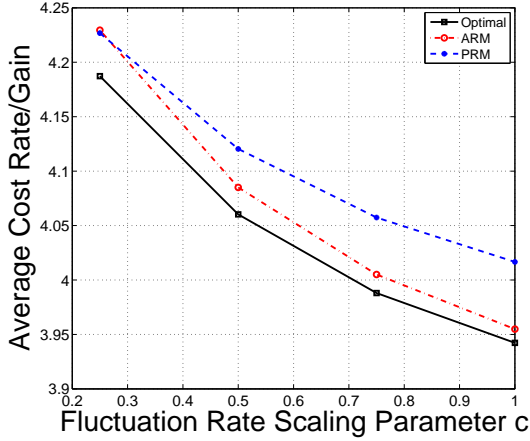}\label{fig:cyc-low}}
\subfigure[Med. Difference in Phase Intensities (Case II)]
{\includegraphics[totalheight=0.28\textheight]{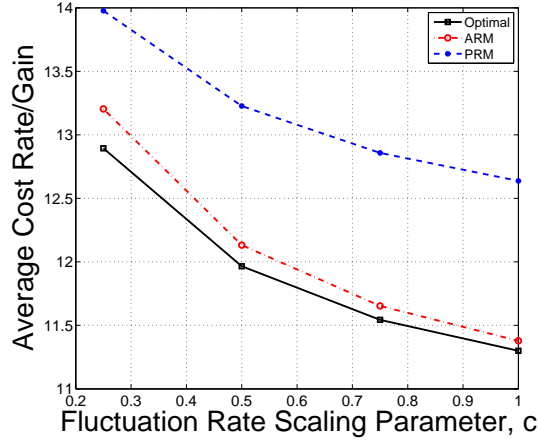}\label{fig:cyc-med}}
\subfigure[Large Difference in Phase Intensities (Case III)]
{\includegraphics[totalheight=0.28\textheight]{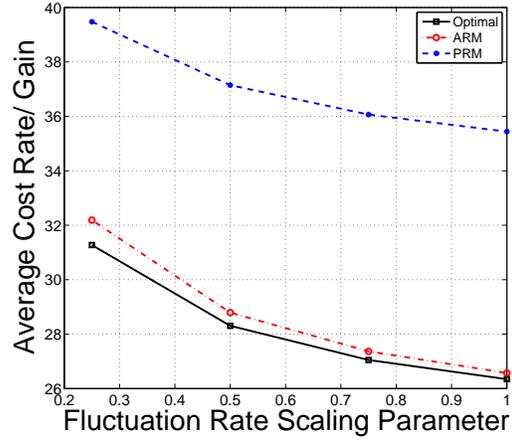}\label{fig:cyc-high}}
\subfigure[Heuristic Vs Optimal Policy (Case III, $n = 2$)]
{\includegraphics[totalheight=0.28\textheight]{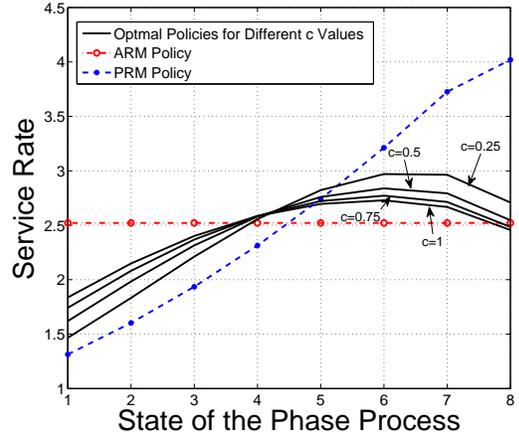}\label{fig:heur-v-opt-cyc}}
\caption{Comparison of Gain Values for Heuristic and Optimal Policies for
Example \ref{ex:numerics2}}
\label{fig:ex_42}
\end{figure}
Figures \ref{fig:cyc-low}-\ref{fig:cyc-high} show the change in the average cost
computed under the heuristics as well as the optimal policies as a function of
parameter $c$ for the three arrival cases. It is interesting to observe that unlike
Example \ref{ex:numerics1}, the ARM policy outperforms the PRM policy in almost
all cases. As illustrated in Example \ref{ex:non-mono}, when the phase process has
cyclic transitions, the optimal service rates may not be monotone in the phase
process (for each fixed $n$). In fact, while the service rates for the PRM policy are
monotone in the phase of the transition process for each congestion level, the
service rates in the optimal policy may begin to decrease as the phase state
increases. This can be more clearly observed in Figure \ref{fig:heur-v-opt-cyc} which
shows a comparison of heuristic and optimal policies as a function of the phase
state when $n=2$ for various values of parameter $c$ and the arrival rates of Case
III. Thus the ARM policy approximates the optimal policy better than the PRM policy
which explains the observed performance difference.

Similar to the previous example, we observe that the performance of the ARM policy
gets worse and that of PRM policy improves with a decrease in the fluctuation rate
parameter $c$. Furthermore, the average cost percent differences provided in Table
\ref{tab:ex_42} show that the ARM policy performs extremely well for all three
arrival rate cases (less than 5\% from optimal for all cases). The PRM policy
performs well for Case I but the degradation in its performance is quite significant
(almost 35 \% from optimal for Case III with $c=1$) when the difference in arrival
rates is medium or high (Cases II and III).

\end{ex}

 Tables \ref{tab:ex_41} and \ref{tab:ex_42} show the optimal costs achievable under the fixed rate mechanism under the column heading ``fixed rate''. We find that costs are between 13\% and 76\% suboptimal when using a fixed rate mechanism relative to the variable rate mechanism. This shows that there is substantial benefit in investing in a responsive mechanism. Furthermore, these examples show that one cannot rely on a particular heuristic method to perform well in all scenarios. In particular, we find that within the gamut of simple heuristic methods considered here, the transition structure and the transition rates of the phase process play an important role in the selection of an appropriate approximation method.

\subsection{Approximation of a System with Non-homogeneous Poisson Arrivals}\label{sec:nhpp}
When the arrival process follows known rate changes, a non-homogeneous Poisson
process is a reasonable modeling tool. In the classical work of Green and Kolesar
\cite{Green1995} or Massey and Whitt \cite{Massey1997} the \textit{analysis} of
queues with non-stationary arrivals is considered. From the standpoint of
\textit{control}, Yoon and Lewis \cite{Yoon2004} consider the case of admission
and pricing control. One thing is certain from Yoon and Lewis's work, control of
non-stationary processes can be computationally intensive. This is due the fact that
to solve each instance the numerical approach requires the time to be discretized.

In this section we explore the possibility of computing an approximate average cost
optimal policy for a single server queue with non-homogeneous Poisson arrivals using
the optimal policies for a system with a ``suitable'' Markov-modulated Poisson arrival
process. Apart from the arrival process, other details are the same as the setting
described in Section \ref{sec:formulate}. Let the arrival process be an NHPP with rate $\lambda(t)$. Assume that $\lambda(t)$ is a periodic function with period $T$. Since the
optimization criterion considered in this study is over an infinite time horizon, and
the rate function for NHPP is a periodic function of time, the \textit{principle of
optimality} implies that only the time elapsed in the current period and the number
of jobs in the system need to be included in the state space \cite{Yoon2004}.

To compute the optimal policy for an NHPP arrival process numerically, the time period is divided into $n$ equally
spaced segments of length $\Delta t=T/n$. Denote the state space for this discretized
process as $\X=\{(n,z)\mid n\in\Z^+,s\in \{0,\Delta t,\dots, T-\Delta t\}\}$. Under
this setting, the decision epochs are the time points $0,\Delta t,2\Delta
t,\dots,T-\Delta t$. Let $\nu$ be a uniformizing rate of the process. An event (arrival,
departure or dummy transition) occurs at a decision epoch with probability
$1-e^{-\nu \Delta t}$ and with probability $e^{-\nu \Delta t}$ no event
occurs. 
The standard theory of Markov decision processes yields the following average cost
optimality inequality (ACOI):
\begin{align*}
	w(n,z)&\geq\min_{x\in\A}\bigg\{\left(-g+h(n)+\mathbf{1}_{\{n>0\}}c(x)\right)\Delta t+
	(1-e^{-\nu \Delta t})\bigg(\frac{\lambda(z)}{\nu} w(n+1,z+\Delta t)\\\nonumber
	&\qquad{}+\mathbf{1}_{\{n>0\}}\frac{x}{\nu}w(n-1,z+\Delta t)+\big(1-\frac{\lambda(z)}{\nu}-\mathbf{1}_{\{n>0\}}\frac{x}{\nu}\big)w(n,z+\Delta t)\bigg)\\\nonumber
	&\qquad{}+e^{-\nu \Delta t}w(n,z+\Delta t)\bigg\}\text{  for 	}n\in\Z^+,z=0,\Delta t,\dots,T-2\Delta t,  \text{ and}\\\nonumber
	w(n,T-\Delta t)\geq&\min_{x\in\A}\bigg\{\left(-g+h(n)+\mathbf{1}_{\{n>0\}}c(x)\right)\Delta t+
	(1-e^{-\nu \Delta t})\bigg(\frac{\lambda(T-\Delta t)}{\nu} w(n+1,0)\\\nonumber
	 &\qquad{}+\mathbf{1}_{\{n>0\}}\frac{x}{\nu}w(n-1,0)+\big(1-\frac{\lambda(z)}{\nu}-\mathbf{1}_{\{n>0\}}\frac{x}{\nu}\big)w(n,0)\bigg)\\\nonumber
	&\qquad{}+e^{-\nu \Delta t}w(n,0)\bigg\}\text{  for 	}n\in\Z^+,\\\nonumber
\end{align*}
where $\mathbf{1}_{E}$ is the indicator function of the event $E$. When the
solution, $(w,g)$ to the ACOI exists, $w$ is called the \textit{relative value
function} and $g^\ast(x)=g$ is the optimal long-run expected average cost for any
initial state $x$.

We now present a method for constructing an approximate policy for NHPP arrivals. The main idea is to approximate the NHPP by an appropriately constructed MMPP. This is done by dividing the time period $T$ into $l$ subintervals and constructing an MMPP with the same number of phases as the number of subintervals i.e, $l$. We choose a cyclic transition structure for the phase process with arrival rate in each phase as the average rate over that subinterval. Transition rates for the phase process can be selected such that the mean sojourn time in phase $s$ is the width of the corresponding interval. The optimal policy corresponding to this MMPP can then be applied to the original NHPP arrivals. The detailed procedure to evaluate the approximate policy is described below:

\begin{enumerate}
\item Partition the interval $[0,T]$ into $l$ subintervals, $[t_{s-1},t_{s}]$, $s = 1,\dots,l$ with $t_0 = 0$ and $t_l = T$.
\item Compute average rates over each partition,

$$
\lambda_s = \frac{\int_{t_{s-1}}^{t_s}\lambda(t)dt}{(t_s-t_{s-1})}.
$$

\item Compute transition rates, $\eta_{i,j}$, for the phase process of MMPP,

$$
\eta_{i,j} = \left\{
\begin{array}{lll}
\frac{1}{(t_i-t_{i-1})}, & 1\leq i\leq l-1 ,j = i+1\hbox{;} \\
\frac{1}{(t_l-t_{l-1})}, & i = l, j = 1 \hbox{;}\\
0, &  \hbox{otherwise}
\end{array}
\right\}
$$.

\item Compute the optimal policy corresponding to the MMPP arrival process constructed in previous steps. Denote this policy as $\mu(n,s)$ , $n\in \Z^+, s \in {1,2,\dots,l}$.
\item Construct the approximate policy, $\hat{\mu}$, for the NHPP process as
$$
\hat{\mu}(n,t) = \mu(n,s) \text{  for } t_{s-1}\leq t \leq t_s \text{ , } n = 0,1,2,\dots.
$$
\end{enumerate}

A numerical study comparing the performance of the approximation procedure stated above for various test
cases is provided in Examples \ref{ex:numerics3} and \ref{ex:numerics4}. In all cases, the policies and average cost are computed using value iteration where the queue length is truncated at 50 to keep the size of the state space manageable.

\begin{ex}
\label{ex:numerics3}
In this example, we consider a NHPP with the following (periodic) rate function,
\begin{align*}
\lambda(t)& =\begin{cases}
		0.1 \text{\quad  for\quad} 0\leq t < T/5\\
		2.0 \text{\quad  for\quad} T/5\leq t < 2T/5\\
		4.0 \text{\quad  for\quad} 2T/5\leq t < 3T/5\\
		2.0 \text{\quad  for\quad} 3T/5\leq t < 4T/5,\\
		0.1 \text{\quad  for\quad} 4T/5\leq t < T,\\
\end{cases}
\end{align*}
where $T$ is the time period of the rate function. One can think of this rate function
as a quantized version of a triangular waveform with time period $T$.
The discretization interval $\Delta t$, for solving the problem with NHPP arrivals is selected as $0.05$
units.  The cost rate function is $c(\mu)=e^{\mu}-1$ and holding cost function is $h(n)=n$. Service rates can be selected from a set $\A=[0,10]$. For computing the approximate policy, we partition the interval $[0,T]$ into 5 subintervals of equal length. Thus, the arrival rates for the corresponding MMPP are $\lambda_1=0.1,\lambda_2=2.0, \lambda_3=4.0$,$\lambda_4=2.0$ and $\lambda_5=0.1$ and the associated generator matrix is
\begin{align*}
\mathbf{Q} & = \frac{5}{T}\left[\begin{matrix}-1&1&0&0&0\\0&-1&1&0&0\\0&0&
-1&1&0\\0&0&0&-1&1\\1&0&0&0&-1\end{matrix}\right],
\end{align*}


\begin{figure}[htbp]
\centering
\includegraphics[totalheight=0.4\textheight]{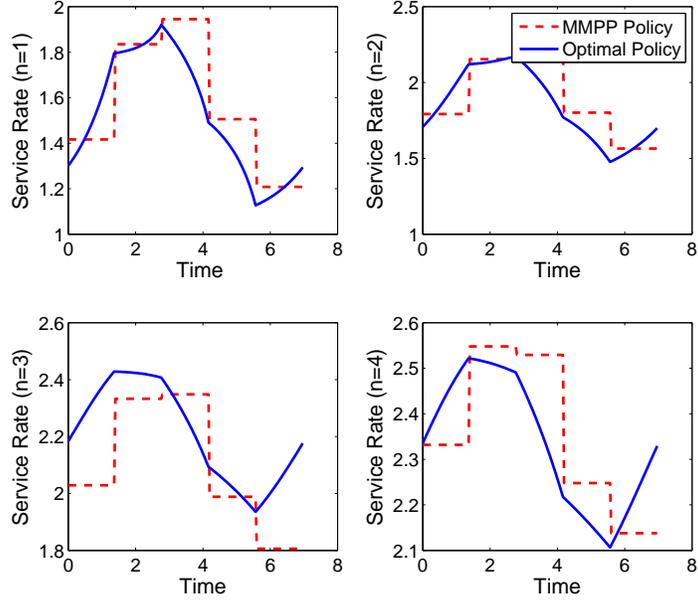}
\caption{NHPP Policy vs Approximate MMPP policy  for various congestion levels (T=7)}
\label{fig:NHPP7}
\end{figure}

\begin{table}[htbp]
\centering
\begin{tabular}{|c|c|c|}
\hline
\multirow{2}{*}{Time Period $(T)$} & \multicolumn{2}{c|}{Gain(\% Sub-Optimal)}  \\
\cline{2-3}
&  Optimal & App. MMPP  \\\hline
4     & 8.5667 & 8.5932 (0.31\%) \\
5     & 8.7750 & 8.7262 (0.41\%) \\
6     & 8.7467 & 8.7925 (0.52\%) \\
7     & 8.8225 & 8.8785 (0.64\%) \\
\hline
\end{tabular}%
\caption{Average Cost Rates and Percentage Difference between Optimal and Approximate NHPP Policy}
\label{tab:NHPP}%
\end{table}%

Figure \ref{fig:NHPP7} show a comparison of optimal policy
with approximate MMPP policy. Table \ref{tab:NHPP} gives the average cost percent
difference between the performance of approximate and optimal policy for various
test scenarios. This data shows that the approximate policies perform extremely well
in all cases (less than 1\% sub-optimal).
\end{ex}

\begin{ex}
\label{ex:numerics4}
In this example, we consider a NHPP with the following (periodic) rate function,
\begin{align*}
\lambda(t)& = 5sin(\omega T) + 6
\end{align*}
where $T$ is the time period of the rate function and $\omega=\frac{2\pi}{T}$ the frequency. We consider the cases in which $T$ is set to $\frac{n\pi}{2}, n = 1,2,3, \text{and}, 4$. The discretization interval $\Delta t$, for solving the problem with NHPP arrivals is selected as $\frac{T}{200}$ units.  The cost rate function is $c(\mu)=\frac{\mu^{2}}{2}-1$ and holding cost function is $h(n)=(n-20)^+$. Service rates can be selected from a set $\A=[0,15]$. For computing the approximate policy, we partition the interval $[0,T]$ into 6 subintervals of equal length.

\begin{figure}[htbp]
\centering
\includegraphics[totalheight=0.4\textheight]{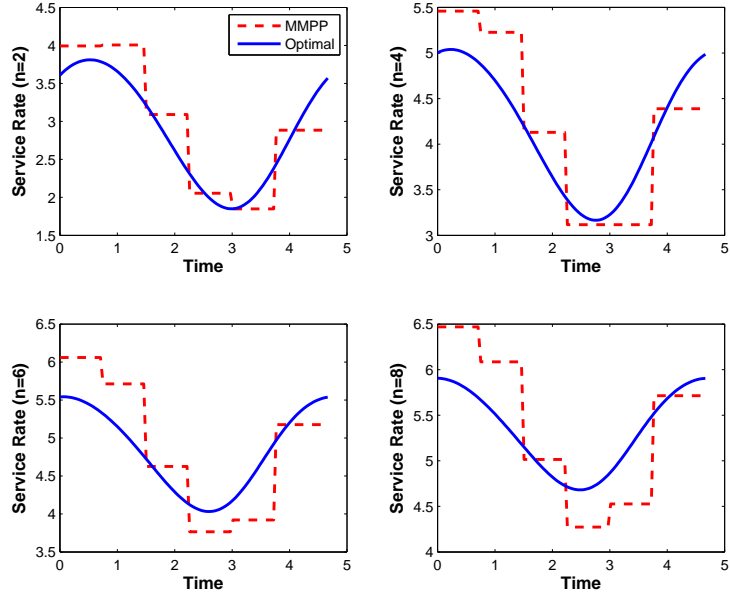}
\caption{NHPP Policy vs Approximate MMPP policy for various congestion levels $(T=\frac{3\pi}{2})$}
\label{fig:NHPPsin}
\end{figure}

\begin{table}[htbp]
\centering
\begin{tabular}{|c|c|c|}
\hline
\multirow{2}{*}{Time Period $(T)$} & \multicolumn{2}{c|}{Gain(\% Sub-Optimal)}  \\
\cline{2-3}
&  Optimal & App. MMPP  \\\hline
$\frac{\pi}{2}$     & 18.90 & 19.38 (2.52\%) \\
$\pi$      & 18.91 & 19.32 (2.21\%) \\
$\frac{3\pi}{2}$     &18.93 & 19.08 (0.82\%) \\
$2\pi$      & 18.94 & 19.29 (1.84\%) \\
\hline
\end{tabular}%
\caption{Average Cost Rates and Percentage Difference between Optimal and Approximate NHPP Policy}
\label{tab:NHPPsin}%
\end{table}%

Figure \ref{fig:NHPPsin} shows a comparison of the optimal policy with heuristic MMPP policy. Table \ref{tab:NHPPsin} gives the average cost percent difference between the performance of approximate and optimal policy for various
test scenarios. We can again see that the approximate policies perform well
in all cases (less than 3\% sub-optimal).
\end{ex}

Of course, we can not conclude from this
limited analysis that MMPP approximate policies will perform well for NHPP in more
general settings. Our motivation in presenting this analysis is to
stimulate further research in this direction. 

\section{Conclusions}\label{sec:conclude}
In this paper we investigate the problem of service rate control for a single server queue with non-stationary arrivals. We propose a framework based on the Markov modulated Poisson processes; a popular model amongst practitioners that is also relatively easy to analyze. Assuming that the goal is to minimize a combination of effort cost and holding cost incurred per unit time, we study this problem under both the discounted and average cost optimality criterion. In either case, we characterize the structure of an optimal service rate as being monotone in the queue length for each arrival rate but \textbf{not} necessarily monotone in the arrival rates for each queue length. In particular, we show that the manner in which the process switches between the arrival rates plays an important role in determining the structure of the optimal policy. We further prove that monotonicity in the arrival rates is recovered when the transition matrix governing the MMPP is stochastically monotone.


There are several ways to extend our work. Our numerical study confirms that in some cases simple heuristics may perform well in the face of changing arrival rates. However, it also points out that careful selection based on the parameters of the system is required, and in many cases applying the proposed model is essential (as opposed to the heuristics). The second part of our numerical work points to the fact that our model can be used as a heuristic itself. We show that we can potentially provide a policy for a system with non-homogenous Poisson arrivals using the optimal policy for an MMPP/M/1 queue. Our results indicate that this may be a promising direction for future research.

Another problem of interest is that of the control of MMPP/M/1 queue with a finite buffer but with an explicit constraint on the job loss rate. Under this setting, while the technical conditions required for stability are not needed, handling the explicit constraint poses a significant challenge. We note that results provided by Ata \cite{Ata2005} for the stationary arrival case may be useful.

The model under study assumes that complete information about arrival statistics is available. This may be unreasonable in situations where arrival statistics cannot be associated with the observable features of the system. A promising direction of work may be to tackle such situations using the partially observable Markov decision process framework.

Finally, we would like to point out that although for ease of exposition, we assume that the cost of effort
function, $c(\mu)$, is strictly convex, continuously differentiable and non-decreasing, our proofs (with minor modification) and results hold for more general cost of effort functions. Using the analysis presented by George and
Harrison \cite{George}, it can be easily shown that the structural results for an optimal policy continue to hold when $c(\mu)$ is assumed to be non-decreasing and continuous. 

\bibliography{Dyn_Cont}
\bibliographystyle{plain}

\section{Appendix}
This appendix is dedicated to providing proofs of Propositions \ref{prop:stable} and
\ref{prop:sennott}. The results of Dai \cite{Dai} and Dai and Meyn \cite{DaiMeyn}
are utilized to show the stability of a stochastic model by establishing the stability of
its fluid limit approximation. For the purpose of this analysis, consider the continuous
time Markov process, $X^\pi(t)=\{(Q^\pi(t),S^\pi(t)),t\geq0\}$, induced by an
admissible stationary policy $\pi\in\Pi$ where $\{Q^\pi(t),t\geq0\}$ and
$\{S^\pi(t),t\geq0\}$ represent the queue length and phase transition process for
arrivals, respectively.

The proof approach for Proposition \ref{prop:stable} follows closely that of Kaufman
and Lewis \cite{Kaufman} (Proposition 3.1) and is done in several steps. Let $\hat{\pi}$ be a policy
that selects a constant rate $\hat{\mu}\in(\sum_{s=1}^L p_s\lambda_s,\bar{u}]$
whenever the queue is non-empty. Let $X^{\hat{\pi}}(t)=\{(Q^{\hat{\pi}}(t),S^{\hat{\pi}}(t)),t\geq0\}$ be the Markov
process induced by the policy $\hat{\pi}$ on state space $\X=\{(n,s)\mid n\in \Z^+,
s\in\mathcal{S}\}$. Since the policy is fixed for the remainder of this section, in the
interest of brevity we suppress dependence on $\hat{\pi}$. The norm of a state,
$x=(n,s)\in \X$, is defined to be $\abs{x} := n+s$. For an initial state $X(0)=x$, we
define the scaled queue length process
\begin{align*}
\bar{Q}^x(t) & := \frac{1}{\abs{x}}Q^x(\abs{x}t).
\end{align*}
We will use a similar notation to denote scaled versions of other stochastic processes.

For each $s\in\mathcal{S}$, let $\{\xi_s(k),k \geq 1\}$ be a sequence of i.i.d
exponential random variables with mean $1/\lambda_s$. The sequence
$\{\xi_s(k),k\geq 1\}$ represents the set of job inter-arrival times when the arrival
process is in phase $s$. Also, let $\{\eta(k),k \geq 1 \}$ be a sequence of i.i.d
exponential random variables with mean $1$ representing the set of job completion
times. Based on these sequences, we define the following cumulative processes
\begin{align*}
		E_s(t) &= \max\{k\geq0\suchthat
\xi_s(1)+\xi_s(2),\dots,\xi_s(k)\leq t\} \qquad \text{for $s\in\mathcal{S}$},\\
		D(t)    &= \max\{k\geq0\suchthat \eta(1)+\eta(2),\dots,\eta(k)\leq  t \}.
\end{align*}
Let $Y^x_s(t)$ be the cumulative amount of time the arrival process spends in phase
$s$ until time $t$ when the initial state is $x$. Similarly, let $T^x(t)$ be the
cumulative amount of time for which there is at least one customer in the queue,
$I^x(t)$ be the cumulative amount of time when the queue is empty and $W^x(t)$
be the total work done by the server until time t. We can now write the following
system of equations for the stochastic process induced by policy $\hat{\pi}$ when
starting from initial state $x$,
\begin{align}
& Q^x(t) =Q^x(0)+\sum_{s=1}^L E_s(Y^x_s(t))-D(W^x(t)), \label{eq:queue}\\
& Q^x(t)\geq0,\\
& \sum_{s=1}^L Y_s^x(t)=t,\label{eq:phase}\\
& W^x(t)=\hat{\mu}T^x(t),\label{eq:work}\\
& T^x(t)+I^x(t)=t,\label{eq:time}\\
& \int_0^\infty Q^x(t)dI^x(t)=0, \label{eq:nonidling} \\
& Y^x_s(t),T^x(t),I^x(t),W^x(t) \quad \text{start from zero and are non-decreasing in $t$.}\label{eq:nondec}
\end{align}
A few comments are in order. First, note that \eqref{eq:nonidling} imposes the
constraint that the server is idle only when the system is empty. For a subsequence
$\{x_n, n \geq 1\}$ such that $|x_n|\rightarrow \infty$, any limit point,
$\bar{Q}(t)$, of the sequence $\{\bar{Q}^{x_n}, n \geq 1 \}$ is called a
\textit{fluid limit}. It will be shown that every fluid limit satisfies a set of equations
known as the \textit{fluid model}. A fluid model is called \textit{stable} if there exists
a $t_0>0$ such that $\bar{Q}(t)=0$ for all $t\geq t_0$ and for all fluid limits. We
next present the fluid model and convergence results for the scaled processes

\begin{prop}
Let $\{x_j\suchthat x_j\in \X, j \geq 1\}$ be a sequence of initial states with
$\abs{x_j}\rightarrow \infty$. Then with probability 1, there exists a subsequence,
$\{x_{j_k}, k \geq 1\}$, such that
	\begin{align}
		&(\bar{Q}^{x_{j_k}}(0),\bar{S}^{x_{j_k}}(0))\rightarrow (\bar{Q}(0),0),\\
    		&(\bar{Q}^{x_{j_k}}(t),\bar{T}^{x_{j_k}}(t))\rightarrow (\bar{Q}(t),\bar{T}(t)) \quad \text{uniformly on compact sets (u.o.c.)},
	\end{align}
	where $(\bar{Q}(t),\bar{T}(t))$ satisfy the following equations,
	\begin{align}
		&\bar{Q}(t)=\bar{Q}(0)+\sum_{s=1}^L p_s\lambda_s t-\bar{W}(t), \label{eq:fmqueue}\\
		&\bar{Q}(t)\geq0,\\
		&\bar{W}(t)=\hat{\mu}\bar{T}(t), \label{eq:fmwork}\\
		&\bar{T}(t)+\bar{I}(t)=t,\label{eq:fmtime}\\
		&\int_0^\infty \bar{Q}(t)d\bar{I}(t)=0, \label{eq:fmnonidling}\\
    		&\bar{T}(t),\bar{I}(t),\bar{W}(t) \quad \text{start from zero and are non-decreasing in $t$}.\label{eq:fmnondec}
	\end{align}
\end{prop}
\Pf  Since $\bar{Q}^{x_j}(0)\leq 1$, $\bar{S}^{x_j}(0)\leq1$ and $1\leq
S^{x_j}(0)\leq L$ for all $j\in\N$, there exists a subsequence, $\{x_{j_k}, k \geq
1\}$ such that $(\bar{Q}^{x_{j_k}}(0),\bar{S}^{x_{j_k}}(0))\rightarrow
(\bar{Q}(0),0)$. For any fixed sample path $\omega$ and $0\leq s \leq t$, we have
$0\leq \bar{T}^{x_j}(t)-\bar{T}^{x_j}(s)\leq t-s$. Thus, the function
$\bar{T}^{x_j}(t)$ is uniformly Lipschitz of order 1. Since $0\leq
\bar{T}^{x_j}(t)\leq t$, it is also uniformly bounded for each $j \geq 1$. Therefore,
the sequence $\{\bar{T}^{x_j}(t), j \geq 1\}$ is equicontinuous. By Arzela-Ascoli
theorem, any subsequence of $\{\bar{T}^{x_j}(t), j \geq 1\}$ has a u.o.c.
convergent subsequence.

Since the phase transition process is ergodic, for each $s \in \{1,\dots,L\}$ we have
with probability 1 $\lim_{t\rightarrow \infty}Y_s^x(t)/t = p_s$.	 Furthermore, from
the strong law of large numbers for renewal processes, the following hold almost
surely
\begin{align*}
	 \lim_{t\rightarrow \infty}E_s(t)/t &= \lambda_s \quad s \in \mathcal{S}, \\
	 \lim_{t\rightarrow \infty}D(t)/t &= 1.
\end{align*}
The above results can be used in a manner similar to Lemma 4.2 of \cite{Dai}, to yield
(with probability 1)
\begin{align}
	\bar{Y}_s(t)& = \lim_{k \to \infty} \frac{1}{\abs{x_{j_k}}}Y^{x_{j_k}}(\abs{x_{j_k}}t)
=p_s t \quad u.o.c.,\text{ for  $s \in \mathcal{S}$}, \label{eq:ssln1}\\
	\bar{E}_s(t)& = \lim_{k \to \infty} \frac{1}{\abs{x_{j_k}}} E(\abs{x_{j_k}}t)
=\lambda_s t \quad u.o.c.,\text{ for  $s \in \mathcal{S}$} \label{eq:ssln2}\\
	\bar{D}(t)& = \lim_{k \to \infty} \frac{1}{\abs{x_{j_k}}}D(\abs{x_{j_k}}t)
=t \quad u.o.c.\label{eq:ssln3}
\end{align}
The equality in \eqref{eq:fmqueue} follows from \eqref{eq:queue} and
\eqref{eq:ssln1}-\eqref{eq:ssln3}. Similarly,
\eqref{eq:fmwork}-\eqref{eq:fmnondec} follow directly from
\eqref{eq:work}-\eqref{eq:nondec}, respectively. \EndPf

\noindent \textbf{Proof of Proposition \ref{prop:stable}:} We start by showing that
the fluid model provided in \eqref{eq:fmqueue}-\eqref{eq:fmnondec} is stable. First
note that $\bar{T}(t),\bar{I}(t) \text{ and } \bar{Y}_s(t)$ are Lipschitz continuous
and therefore absolutely continuous and differentiable almost everywhere. Taking the
derivative with respect to $t$ in \eqref{eq:fmqueue} and \eqref{eq:fmwork} yields
\begin{align*}
		\dot{\bar{Q}}(t) & =\sum_{s=1}^L p_s\lambda_s-\hat{\mu}\dot{\bar{T}}(t).
\end{align*}
Further, due to the non-idling constraint \eqref{eq:fmnonidling},
$\dot{\bar{I}}(t)=0$ whenever $\bar{Q}(t)>0$. Thus, from \eqref{eq:fmtime},
$\dot{\bar{T}}(t)=1$ whenever $\bar{Q}(t)>0$. So for $\bar{Q}(t)>0$, we have
\begin{align*}
	\dot{\bar{Q}}(t) & =\sum_{s=1}^L p_s\lambda_s-\hat{\mu}.
\end{align*}
Our choice of the stationary policy enabled by the stability condition
\eqref{eq:stability}, implies that $\dot{\bar{Q}}(t)<0$ whenever $\bar{Q}(t)>0$.
Thus from Lemma 5.2 of \cite{Dai} we have that the fluid limit process for queue
length is non-increasing and there exists a $t_0\geq 0$ such that $\bar{Q}(t) = 0$
for all $t\geq t_0$. That is, the fluid model is stable. The results of Theorem 4.2 of
\cite{Dai} imply that the Markov process induced by the stationary policy
$\hat{\pi}$, is positive recurrent and stationary distribution exists. Furthermore,
since the embedded discrete time Markov chain for the process is irreducible, this
process is \textit{ergodic}. This implies that the long-run average cost under
$\hat{\pi}$, say $g^{\hat{\pi}}$, is independent of the initial state.

To show that $g^{\hat{\pi}}$ is finite, we use the results of Theorem 4.1(i) of
\cite{DaiMeyn}. Since the fluid model is stable and conditions A1) and A2) of
Theorem 4.1 hold, it follows that for any integer $p\geq 1$
$\lim\sup_{t\rightarrow\infty}\frac{1}{t}\int_0^t\E_x[|Q(u)|^p]du<\infty$ for each
initial condition $x$. Since the holding cost has a polynomial rate of growth, we have that the long-run average
holding cost rate is finite. Moreover, the direct contribution to long-run cost rate due
to serving at $\hat{\mu}$ whenever the queue is not empty is at most
$c(\hat{\mu})<\infty$. It therefore follows $g^{\hat{\pi}}$ is finite. 	

It remains to consider the necessity of \eqref{eq:stability}. Consider the Markov
process induced by a policy that uses the highest available service rate whenever the
queue is not empty. As shown by Yechiali\cite{Yechiali} using the detailed balance
equations for the steady state distribution a non-trivial invariant measure exists only
if $\bar{u}>\sum_{s=1}^L p_s\lambda_s$. The result follows and the proof is
complete. \EndPf

The remainder of this section is dedicated to proving Proposition \ref{prop:sennott}
holds. We show that the optimal value and relative value functions satisfying the
ACOI, \eqref{eq:ACOI}, exist and can be obtained via limits from the discounted
expected cost value functions. Thus, the structural results proved for the discounted
cost case continue to hold for the average cost case. In proving these results, we
verify the following set of assumptions $(SEN)$ (included for completeness) provided
by Sennott \cite{sennott89} hold.
\begin{itemize}
	\item \textit{SEN1:} There exist $\delta>0$ and $\epsilon>0$ such that for every state and action, there is a probability of at least $\epsilon$ that the
transition time will be greater that $\delta$.
	\item \textit{SEN2:} There exists $B$ such that $\tau(i,\mu)\leq B$ for every state $i$ and control $\mu$, where $\tau(i,\mu)$ is the expected transition time out of state $i$ when control $\mu$ is chosen.
	\item \textit{SEN3:} $v_\alpha(i)<\infty$ for all states $i$ and $\alpha>0$.
	\item \textit{SEN4:} There exists $\alpha_0>0$ and nonnegative numbers $M_i$ such that $w_{\alpha}(i)\leq M_i$ for every state $i$ and $0<\alpha<\alpha_0$ where $w_\alpha(i)=v_\alpha(i)-v_\alpha(\mathbf{0})$, for distinguished state $\mathbf{0}$. For every state $i$, there exists an action $\mu_i$ such that $\sum_jP_{ij}(\mu_i)M_j<\infty$.
	\item \textit{SEN5:} There exists $\alpha_0>0$ and a non-negative number $N$ such that $-N\leq w_{\alpha}(i)$ for every $i$ and $0\leq\alpha\leq \alpha_0$.
	\item \textit{SEN6:} For each state $i$, the expected single stage discounted
cost $f_\alpha(i,\mu)$ is a lower semi-continuous (lsc) function on the product
space $[0,\infty)\times\A$. Note that $f_0(i,\mu)$ is the un-discounted single
stage cost.
	\item \textit{SEN7:} For all states $i$ and $j$, the function
$L_{ij}(\alpha,\mu)=P_{ij}(\mu)\int_{t=0}^{\infty}e^{-\alpha t}\nu e^{-\nu
t}dt$ is a  lsc function on the product space  $[0,\infty)\times\A$.
	\item \textit{SEN8:} Assume that $\alpha_n$ is a sequence of discount factors
converging to 0 with the property that $\pi_{\alpha_n}$, the associated sequence
of $\alpha$-discount optimal policies, converge to a stationary policy $\pi$. Then
for each state $i$, $\lim\inf_{n}\tau(i,\pi_{\alpha_n})\leq \tau(i,\pi)$.
\end{itemize}

\noindent \textbf{Proof of Proposition \ref{prop:sennott}:} We begin by verifying
the \textit{SEN} assumptions. Since the uniformizing rate is strictly positive and
finite, assumptions \textit{SEN1} and \textit{SEN2} hold. It was shown in
Proposition \ref{prop:stable} that under the stability condition \eqref{eq:stability},
there exists a stationary policy that induces an ergodic Markov process with finite
long-run average expected cost. Thus, the hypotheses of Lemma 2 of
\cite{sennott89} hold; validating \textit{SEN3 and SEN4} assumptions. Let
$\bar{s}=\arg\min_{s\in\mathcal{S}}\{v_\alpha(0,s)\}$ and define the
distinguished state as $\mathbf{0}=(0,\bar{s})$. It follows from Proposition 3.1
that for any $\alpha>0$, $v_\alpha(\mathbf{0})\leq v_\alpha(n,s)$ for all
$(n,s)\in\X$. Therefore, $w_\alpha(n,s)\geq0$ and \textit{SEN5} is satisfied.
\textit{SEN6} holds since for each $(n,s)\in \X$, $f_{\alpha}(
(n,s),\mu)=(c(\mu)+h(n))/(\alpha+\nu)$ is a continuous function on
$[0,\infty)\times\A$.

There is no decision to be made when $n=0$. Fix $n\geq1$ and $s\in\{1,2,\dots,L\}$
and note,
\begin{align}
	L_{(n,s),(n^\prime,s^\prime)}(\alpha,\mu)& =
	\begin{cases}
		\frac{\lambda_s}{\alpha+\nu}  \text{\quad if\quad } n^\prime = n+1, s^\prime=s,\\
		\frac{\mu}{\alpha+\nu}  \text{\quad if\quad}     n^\prime = n-1, s^\prime=s,\\
		\frac{Q_{ss^\prime}}{\alpha+\nu}  \text{\quad if\quad}     n^\prime = n, s^\prime\in\mathcal{S},\\
		0 \text{\qquad otherwise.\quad}
	\end{cases}
\end{align}
So the function $L_{x,x^\prime}(\alpha,\mu)$ is jointly continuous in $\alpha$ and
$\mu$ for each $x,x^\prime\in \X$, therefore \textit{SEN7} holds. Finally, since for
any policy $\pi$,  $\tau(i,\pi)=1/\nu$, \textit{SEN8} holds. \par
\textit{SEN1}, \textit{SEN3}, \textit{SEN6} and \textit{SEN7} are required in
Theorem 11 of \cite{sennott89} to prove the first result. Since the holding costs are
assumed to have a polynomial rate of growth, the hypotheses of Proposition 4 of \cite{sennott89} are
satisfied. The last two results follow from Theorem 12 (and its proof) of
\cite{sennott89}. \EndPf

\end{document}